[1994/12/01]% LaTeX date must December 1994 or later
\documentclass[reqno,13pt]{amsart}
\usepackage{amsmath,amssymb,amsfonts}

\usepackage{bbm}
\def\1{\mathbbm{1}}

\usepackage[mathscr]{euscript}
\usepackage{graphicx}
\usepackage{pstricks}
\usepackage{pst-node}

\usepackage{tikz}
\usetikzlibrary{arrows,chains,matrix,positioning,scopes}
\tikzstyle{element}=[rectangle,draw,fill=white, line width=1pt]
\tikzstyle{terminal}=[circle,draw, scale=0.3, line width=1pt,red]
\tikzstyle{fleche}=[->,>=stealth', very thick]
\tikzstyle{fleche1}=[->,>=stealth', very thick, red]
\usepackage{placeins}

\def\R{{\mathbb{R}}}

\def\t{\tau}

\def\B{\mathbb{B}}

\def\calM{{\mathcal{M}}}
\def\calR{{\mathcal{R}}}
\def\calL{{\mathcal{L}}}

\def\calT{{\mathcal{T}}}

\def\calA{{\mathcal{A}}}

\def\calC{{\mathcal{C}}}

\def\R{\mathbb R}
\def\calN{\mathcal{N}}

\def\C{\mathbb C}

\def\N{\mathbb N}

\def\F{\mathbb F}

\def\T{\mathbb T}

\newtheorem{definition}{Definition}[section]

\newtheorem{lemma}{Lemma}[section]
\newtheorem{proposition}{Proposition}[section]
\newtheorem{corollary}{Corollary}[section]
\newtheorem{theorem}{Theorem}[section]
\theoremstyle{remark}
\newtheorem{remark}{Remark}[section]
\newtheorem{example}{Example}[section]

\title{Boundary approximate controllability under positivity constraints of linear systems  }
\author{Yassine El Gantouh} 
\address{Departamento de Matem\'{a}ticas, Universidad Aut\'{o}noma de Madrid, 28049 Madrid, Spain.}
\address{Chair of Computational Mathematics, Fundación Deusto, Universityof Deusto, 48007 Bilbao, Basque Country, Spain.}
\thanks{The author wish to thank Prof. Enrique Zuazua for his comments and suggestions.}
\thanks{This work has been supported by the European Research Council (ERC) under the European Union’s Horizon 2020 research and innovation programme (grant agreement No 694126-DYCON).}% At most 5 thanks

\begin{document}

	%	\subtitle{Boundary approximate controllability under positivity constraints}

	\maketitle
	
		\renewcommand{\sectionmark}[1]{}
	
	\begin{abstract}
		This paper focuses on boundary approximate controllability under positivity constraints of a wide range of infinite-dimensional control systems. We develop frequency domain controllability criteria. Firstly, we derive a controllability result under positivity constraints on the control for such systems. Then, and more importantly, we provide a necessary and sufficient condition for controllability under positivity constraints on the control and the state.
		
		The obtained results are applied to the controllability of transportation and heat conduction networks. In particular, provided that the underlying graph is strongly connected, the controllability under positivity constraints on the control/state of transport network systems is fully characterized by a Kalman-type rank condition. For a system of heat equations with Robin boundary conditions on a path-like network, we establish approximate controllability under positivity state-constraint with a single positive input through the starting node. However, we prove the lack of controllability under unilateral control-constraint.
	\end{abstract}
%	\keywords{Infinite-dimensional control systems \and  positive semigroups \and controllability under constraints \and  PDEs on networks}
%	\subclass{93C20 \and  93B05 \and 35R02}
	
	%All acknowledgements should be placed in the back of the paper after Conclusions..
	
	\section{Introduction}\label{Sec0}
	
	In this paper, we are concerned with boundary approximate controllability under positivity constraints of infinite-dimensional linear systems described as
\begin{equation}\label{S5.Sy1}
	\left\lbrace
	\begin{array}{lll}
		\dot{z}(t)= A_{m}z(t), \qquad \qquad\qquad t> 0,\\
		z(0)=x,\\
		(G-\Gamma)z(t)=Ku(t), \qquad \quad  t> 0,
	\end{array}
	\right.
\end{equation}
where the state variable $ z(.) $ takes values in a Banach space $ X $ and the control function $u(\cdot)$ is given in the Banach space $ L^{p}([0,+\infty);U) $. The maximal (differential) operator $A_m:D(A_m)\subset X\to X$ is closed and densely defined, $ K $ is a bounded linear operator from $  U$ into $ \partial X $ (both are Banach space), and $ G,\Gamma: D(A_m)\to \partial X$ are linear continuous trace operators. Such a system arises naturally as abstract formulation of structured population models, and systems of linear partial differential and/or delay differential equations on networks. In particular, over the past few decades, there has been a growing interest in studying qualitative and quantitative properties of (\ref{S5.Sy1}) (see e.g. \cite{BBEM,ChM,CZ,EBMB,El,HMR,LT,Lions,Salam,Sta} and references therein). Before going further in our exposition, let us recall some basic facts about this kind of systems and their controllability properties:
\begin{itemize}
	\item[1.] System (\ref{S5.Sy1}) is well-posed in the following sense: for every $x\in X$ and $u\in L^p([0,+\infty);U)$ there exists a unique solution $z\in C([0,+\infty);X)$ that depends continuously on the initial data $x$ and the control $u$.
	\item[2.] System (\ref{S5.Sy1}) is positively well-posed if it is well-posed and for every positive initial state $x$ and positive input $u$ the state of (\ref{S5.Sy1}) remains positive for all $t\ge 0$. 
	\item[3.] System (\ref{S5.Sy1}) is said to be approximately controllable with respect to positive controls if it can be steered from any initial data to a state arbitrarily close to a target by choosing a suitable positive control. 
	\item[4.] System (\ref{S5.Sy1}) is said to be approximately positive controllable if it can be steered from any positive initial data to a positive state arbitrarily close to a target by choosing a suitable positive control.
\end{itemize}
The well-posedness of (\ref{S5.Sy1}) has been investigated in several works (see, e.g., \cite{ChM,El,HMR,LT,Lions,Salam}). Moreover, boundary approximate positive controllability for the unperturbed system (\ref{S5.Sy1}) (i.e., $\Gamma=0$) has been addressed in \cite{BBEM}. In this latter, the authors proved a necessary and sufficient condition for this system to be approximate positive controllable by developing the solution into an implicit variation of constant formula \cite[Eq. 3.2]{BBEM}. In this setting, it is clearly pointed out that some bounds on the trace operator $G$ lead unavoidably to non-reflexive Banach spaces. More recently in \cite{EBMB}, the authors also addressed the question of boundary approximate controllability with respect to positive controls for a linearized Saint-Venant equation and a heat equation. There, the authors established that the resulting problem can be reduced to a classic constrained controllability problem for distributed systems \cite{Son}. Note, however, that, in this paper, the analyticity of the semigroup were necessary to obtain controllability results. Finally, for the sake of completeness, we mention that a controllability problem under positivity constraints for the heat equation was recently addressed in \cite{LTZ}. There, the authors showed that the heat equation is controllable to any positive steady state by means of positive boundary controls, provided the control time is large enough. Moreover, it is also proved that controllability by positive controls fails if the time is too short, whenever the initial datum differs from the final target.

In this paper, we investigate boundary approximate controllability under positivity constraints of (\ref{S5.Sy1}) within the framework of positive $L^p$-well-posed and regular linear systems introduced in \cite{El}. Our aim is to derive controllability criteria for (\ref{S5.Sy1}) and to generalize and unify some previous results available in the literature \cite{BBEM,Brammer,Klamka,Tilman,Son}. To this end, we derive necessary and sufficient conditions for the boundary approximate controllability  with respect to positive controls of the non-homogeneous boundary control system (\ref{S5.Sy1}). More precisely, we provide frequency domain criteria for boundary approximate controllability under positivity constraints on the control and/or the state. This latter effort is motivated by the fact that transfer functions of infinite-dimensional control systems provide a very clear characterization of the qualitative and quantitative properties of boundary control systems and time-delay systems.  

In particular, we have obtained the following controllability results:
\begin{itemize}
	\item[(i)] Firstly, we show in Theorem \ref{Main-T1} that boundary approximate controllability under positivity constraints on the control is fully characterized by a frequency domain-type test.
	\item[(ii)] Secondly, in Theorem \ref{Main-T2} we establish that boundary approximate controllability under positivity constraints on the control and the state (also called boundary approximate positive controllability) is equivalent to a frequency domain test (inequality-type test).
\end{itemize}

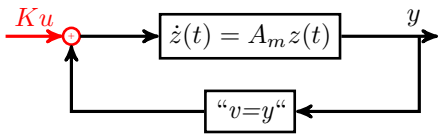
\begin{figure}[ht]\label{Fig1}
	\begin{center}
		\begin{tikzpicture}[scale=0.5]
			\node[terminal] (p) at (-0.75,0) {$\textcolor{red}{\textbf{+}}$};
			\node[element] (S) at (4,0) {$\dot{z}(t)= A_mz(t)$};
			\node[element] (F) at (4,-2) {$``\emph{v=y}``$};
			\draw[fleche] (p) -- (S) node[midway, above]{};
			\draw[fleche] (S)-- (9,0) node[near end, above]{$y$};
			\draw[fleche] (8.5,0) |-(F);
			\draw[fleche] (F) -| (p);
			\draw[fleche1] (-2.5,0) -- (p) node[midway, above]{$\textcolor{red}{Ku}$};
		\end{tikzpicture}
	\end{center}	
	\caption{A closed-loop representation of (\ref{S5.Sy1}).}
	\vspace{.1 in}
\end{figure}

The set-up of the proof of Theorems \ref{Main-T1} and \ref{Main-T2} is relatively simple. First, we assume that the homogeneous system is  positively well-posed (rep. well-posed). Using a perturbation result in \cite{El} (resp. \cite{HMR}), we have that the generator of the homogeneous system coincides with the generator of the closed-loop system associated to (\ref{S5.Sy1}) (see Figure \ref{Fig1}). Therefore, we express the control map of (\ref{S5.Sy1}) in terms of the input-maps of the closed-loop system. This allows us to use Laplace transform arguments (\cite{WR}) and thus obtain frequency domain tests for boundary approximate controllability under positivity constraints. 

The rest of this paper is organized as follows: we state the main results of the paper in Section \ref{Sec1}. In Section \ref{Sec2}, we give a brief exposition of the concept of positive $L^p$-well-posed and regular linear systems. Moreover, we show in this section the existence of positive mild solution of \eqref{S5.Sy1} (Theorem \ref{wellposedness}). Section \ref{Sec3} is devoted to the proof of the main results. The fourth section deals with the application of the results obtained in the previous sections to specific partial differential equations (PDEs) on networks including transport network systems with controls at internal and/or external vertices, and a system of heat equations on a directed network with controls in Robin boundary conditions. For transport networks, if the underlying graph is strongly connected, boundary approximate positive controllability is fully characterized by a Kalman-type rank condition (Corollary \ref{COR1}). In particular, we obtain that transport processes on directed cycles are approximately positive controllable with a control acting on the starting node (Example \ref{EX2}). For a system of heat equations with Robin boundary conditions on a path-like network, we establish approximate controllability under positivity state-constraint with a single positive input through the starting node (Lemma \ref{Lemma-path1}). However, we prove the lack of controllability under unilateral control-constraint (Lemma \ref{Lemma-path}). Finally, we include in the Appendix two technical lemmas needed for proving the Kalman-type rank characterization of boundary approximate positive controllability of transportation networks.

	\section{Main results}\label{Sec1}
	
 Before stating the main results of this paper, we need to introduce some notations and concepts. Let $X$ be a real vector space and $\leq$ be a partial order on this space. Then $X$ is said to be an ordered vector space if it satisfies the following properties:
\begin{itemize}
	\item[(\emph{i})] If $f,g\in X$ and $f\leq g$, then  $f+h\leq g+h$ for all $h\in X$.
	\item[(\emph{ii})] If $f,g\in X$ and $f\leq g$, then  $\alpha f\leq \alpha g$ for all $\alpha\geq 0$.
\end{itemize}
If, in addition, $X$ is lattice with respect to the partial ordered, that is, $\sup\{f,g\}$ and $\inf\{f,g\}$ exist for all $f,g\in X$, then $X$ is said to be a \emph{vector lattice}. For an element $f$ of a \emph{vector lattice} $X$, the \emph{positive part} of $f$ is defined by $f_+:=\sup\{f,0\}$, the \emph{negative part} of $f$ by $f_-:=\sup\{-f,0\}$ and the absolute value of $f$ by $\vert f\vert:=\sup\{f,-f\}$, where $0$ is the zero element of $X$. The set of all positive elements of $X$, denoted by $X_+$, is a convex cone with vertex $0$. In particular, it generates a canonical ordering $ \leq $ on $ X $ which is given by: $f\leq g$ if and only if $g-f\in X_+$. A linear subspace $\mathbb{B}$ of a vector lattice $X$ is said to be \emph{ideal}, if $f\in X,\; g\in \mathbb{B}$, and $\vert f\vert \leq \vert g\vert$ implies $g\in \mathbb{B}$. A norm complete vector lattice $ X$ such that its norm satisfies the following property
\begin{displaymath}
	\vert f\vert \leq \vert g\vert \quad \Longrightarrow \quad \Vert f\Vert \leq \Vert g\Vert,
\end{displaymath}
for $f,g\in X$, is called \emph{Banach lattice}. If $X$ is a Banach lattice, its topological dual $X'$, endowed with the dual norm and the dual order, is also a Banach lattice. We define the polar cone of a subset $ M \subseteq X $ by
$$ 
M^{\circ}=\left\{ \varphi\in X',\;\; \langle \varphi\,,\,f\rangle \leq 0,\;\; \forall f\in M\right\}.
$$
We have the following density result, see \cite[Proposition 2.3]{Son}.
\begin{lemma}\label{Main-lemma}
	Let $X$ be a Banach space and $M$ a convex cone in $X$ with vertex at the origin. Then $\overline{M}=X$ if and only if $M^{\circ}=\{0\}$.
\end{lemma}
We denote by $ \mathcal{L}(E,F) $ the Banach algebra of all linear bounded operators from a Banach space $E$ to a Banach space $F$. If $Z$ is a subspace of a Banach space $E$, then by $P_{\vert Z}$ we denote the restriction operator of $P\in \mathcal{L}(E,F) $. An operator $ P\in \mathcal{L}(E,F)  $ is positive if and only if $ P E_+\subset F_+ $ or, equivalently, if $ f\leq g $ implies $ Pf\leq Pg $. An everywhere defined positive operator from a Banach lattice to a normed vector lattice is bounded, see e.g. \cite[Theorem II.5.3]{Schaf}. The set of all positive operators from a Banach lattice $E$ to another Banach lattice $Y$, denoted by $ \mathcal{L}_+(E,F) $, is a convex cone in the space $ \mathcal{L}(E,F) $. 

Let $X$ be a real Banach lattice and $ (A,D(A)) $ be the generator of a C$_{0} $-semigroup $ \mathbb{T}:=(T(t))_{t\geq 0} $ on $ X $. The type of $ \mathbb{T} $ is defined by $ \omega_{0}(A):=\inf \lbrace t^{-1}\log \Vert T(t)\Vert: t> 0\rbrace $. We denote by $ \rho (A) $ the resolvent set of $ A $, i.e., the set of all $ \mu \in \C $ such that $ \mu I_X -A $ is invertible with $I_X$ denote the identity operator in $X$. By $ R(\mu,A): =(\mu I_X-A)^{-1} $ the resolvent operator of $ A $. The complement of $ \rho (A) $, is called the spectrum and is denoted by $\sigma(A):= \C\backslash \rho (A)$. The so-called spectral radius of $A$ is defined by $r(A):=\sup\{\vert \mu\vert :\; \mu \in \sigma(A)\}$. Also, recall that the spectral bound $s(A)$ of $A$ is defined by $s(A):=\sup\{{\rm Re}\, \mu :\; \mu \in \sigma(A)\}$. A linear operator $A$ on a Banach lattice $X$ is called resolvent positive if there exists $\omega\in \mathbb{R}$ such that $(\omega,\infty)\subseteq \rho (A) $ and $R(\mu,A)\geq 0$ for each $\mu> \omega$. Moreover, by \cite[Corollary 2.3]{Arendt} a C$_{0} $-semigroups on a Banach lattice is positive if and only if the corresponding generator $A$ is resolvent positive. On the other hand, by $ X_{1} $ we denote the order Banach space $ D(A) $ endowed with the norm $ \Vert x\Vert_{1}:=\Vert (\mu I_X-A)x\Vert $ for some $ \mu \in \rho (A) $. The extrapolation space associated with $ X $ and $ A $, denoted by $ X_{-1} $, is the completion of $ X $ with respect to the norm  $\Vert x\Vert_{-1}:=\Vert R(\mu , A)x\Vert $ for $ x\in X $ and some $ \mu \in \rho (A) $. Note that the choice of $ \mu  $ is not important, since by the resolvent equation different choices lead to equivalent norms on $ X_{1} $ and $ X_{-1} $. Moreover, we have
$
X_{1}\subset X\subset X_{-1}.
$
The unique extension of $ \mathbb{T} $ on $ X_{-1} $ is a $C_{0}$-semigroup which we denote by $ \mathbb{T}_{-1}:=(T_{-1}(t))_{t\geq 0} $ and whose generator is denoted by $ A_{-1} $. For more details on positive semigroups, see for instance \cite{BFR} or \cite{Nagel}. 

Let $ X,\partial X, U $ be Banach lattices called \emph{the sate space}, \emph{boundary space} and \emph{input space}, respectively. We rewrite (\ref{S5.Sy1}) as the following boundary input-output system
\begin{equation}\label{input-output}
	\left\lbrace
	\begin{array}{lll}
		\dot{z}(t) =A_{m} z(t),\qquad \qquad \quad \, t> 0,\qquad z(0)=x,\\
		G z(t)-v(t)= Ku(t),\qquad  t> 0,\\
		y(t) = \Gamma z(t), \qquad \qquad \qquad t> 0,
	\end{array}
	\right.
\end{equation}
with the feedback law
\begin{equation}\label{feedback}
	``\emph{v=y}".
\end{equation}
Moreover, let us assume that the trace operator $G$ satisfies:
\begin{itemize}
	\item[{\bf (H1)}]   $ A=A_{m} $ with domain $ D(A):=\ker G $ generates a C$_0$-semigroups $\mathbb{T}$ on $X$;
	\item[{\bf (H2)}]  $ G $ is surjective.
\end{itemize}
According to \cite[Lemma 1.2]{Gr}, the assumptions \textbf{($  \mathbf{H1}$)}-\textbf{($  \mathbf{H2}$)} imply that the domain of $ A_m $ can be decomposed and related to $ A $ as
$$
D(A_m)=D(A) \oplus \ker ( \mu I_X- A_{m} ),\qquad \mu \in \rho(A).
$$
Moreover, the restriction operator $G_{\vert_{ \ker ( \mu I_X-A_{m} )}}$ is invertible and the operator $$ D_{\mu}:= \left(G_{\vert_{ \ker ( \mu I_X-A_{m} )}}\right)^{-1} ,$$ called the Dirichlet operator, exists and is bounded. Thus,
\begin{equation}\label{S5.5}
	A_{m}=A_{-1}+ BG,
\end{equation}
since $\mu D_\mu v=A_m D_\mu v,\,v\in \partial X$, where $B:\partial X\to X_{-1}$ is defined by  
$$
B=(\mu I_X-A_{-1})D_\mu,\qquad \mu\in \rho(A).
$$
Note that $B$ does not depends on $\mu$, due to the resolvent equation.

We are now ready to state the main controllability results of this paper. The first one highlights boundary approximate controllability under positivity control-constraint.
\begin{theorem}\label{Main-T1}
	Let $X,\partial X,U$ be order Banach spaces and let the assumptions {\bf (H1)}--{\bf (H2)} be satisfied. Furthermore, assume that:
	\begin{itemize}
		\item[{\bf (H3)}] $ (A,B,\Gamma_{\vert_{D(A)}}) $ is an $L^p$-well-posed regular triplet (with feedthrough zero) on $ X,\partial X, \partial X $ with the identity $I_{\partial X} $ as an admissible feedback operator.
	\end{itemize}
	Then, the system (\ref{S5.Sy1}) is boundary approximately controllable with respect to positive controls if and only if there exists $\omega\in \R$ such that, for all $\mu\ge \omega$, we have $1\in \rho(\Gamma D_\mu)$ and
	\begin{equation}\label{charactera-T1}
		\left(D_{\mu}(I_{ X}- \Gamma D_{\mu})^{-1}KU_{+}\right)^{\circ}=\{0\}.
	\end{equation}
\end{theorem}

The second and most important result of the present paper provides a complete description of controllability under positivity constraints on the control and the state.
\begin{theorem}\label{Main-T2}
	Let $X,\partial X,U$ be Banach lattices and let the assumptions {\bf (H1)}--{\bf (H2)} be satisfied. Furthermore, assume that $K\in  \calL(U,\partial X)$ is positive and 
	\begin{itemize}
		\item[{\bf (H3')}] $ (A,B,\Gamma_{\vert_{D(A)}}) $ is a positive $L^p$-well-posed regular triplet (with feedthrough zero) on $ X,\partial X, \partial X $ with $ I_{\partial X} $ as a positive admissible feedback operator.
	\end{itemize}
	Then, the system (\ref{S5.Sy1}) is boundary approximately positive controllable if and only if there exits $\omega>s(A)$ such that $ r(\Gamma D_{\omega})<1$ and for all $\mu\ge  \omega$ we have
	\begin{equation}\label{characteraT2}
		\bigcap\left\{\big(D_{\mu}(\Gamma D_{\mu})^n KU_{+}\big)^{\circ},\;\; n\in \N \right\} =X^{\circ}_+.
	\end{equation}
\end{theorem}

\begin{remark}\label{Main-Remark}
	The two theorems above provide frequency-domain characterizations of the boundary approximate controllability under positivity constraints of (\ref{S5.Sy1}). These characterizations are given in terms of transfer functions of the input-output system (\ref{input-output}). On the one hand, Theorem \ref{Main-T1} shows that the perturbed boundary control system (\ref{S5.Sy1}) is boundary approximately controllable with respect to positive controls if and only if there exists $\omega\in \R$ such that, for all $\mu\ge \omega$, $1\in \rho(\Gamma D_\mu)$ and the following implication holds for all $\varphi\in X'$:
	$$
	\langle D_{\mu}(I_{ X}- \Gamma D_{\mu})^{-1}Ku,\varphi\rangle \leq 0, \qquad\, \forall\, u \in U_+,\; \mu\ge \omega\quad \Longrightarrow \quad\varphi=0.
	$$
	On the other hand, Theorem \ref{Main-T2} shows that, under the positivity of the operators $\Gamma,D_\mu,K$, (\ref{S5.Sy1}) is boundary approximately positive controllable if and only if there exits $\omega>s(A)$ such that $r(\Gamma D_{\omega})<1$ and the following implication holds for all $\mu> \omega$ and $\varphi\in X'$:
	$$
	\langle D_{\mu}(\Gamma D_{\mu})^n Ku,\varphi\rangle \leq 0, \qquad\, \forall\, u \in U_+,\,n\in \N,\;\mu\ge \omega\quad \Longrightarrow \quad \varphi \le  0.
	$$
	It is to be noted that if for some $\omega>s(A)$, $r(\Gamma D_{\omega})<1$, then $r(\Gamma D_{\mu})<1$ for all $\mu \ge \omega$ as the family $(\Gamma D_{\mu})_{\mu>s(A)}$ is positive and monotonically decreasing. Finally, let us point out that in this work we do not request any spectral properties on the operators $A$, in contrast to the papers \cite{EBMB,Son} where the existence of a Riesz basis of generalized eigenvectors or related spectral properties are required.
\end{remark}
	
	\section{Well-posedness}\label{Sec2}

  In this section we investigate the well-posedness and positivity property for the solution of \eqref{S5.Sy1}. To this end, we shall use the feedback theory of positive $L^p$-well-posed and regular linear systems developed in \cite{El}. In fact, let $X,\, \partial X,\, U$ be Banach lattices. For $ \alpha>s(A)$, let $ L^{p}_{\alpha}(\R_+;U) $ denote the space of all functions of the form $ v(t)=e^{\alpha t}u(t) $, where $ u\in L^{p}(\R_+;U) $. Moreover, let $L_{loc,+}^p(\R_+;U)$ denote the set of positive control functions $u$ in $L_{loc}^p(\R_+;U)$ such that $u(t)\in U_+$ almost everywhere in $\R_+$, where we regard $L_{loc}^p(\R_+;U)$ as a lattice ordered Fr\'{e}chet space with the seminorms being the $L^p$ norms on the intervals $[0,n]$, $n\in\N$.

We select the following definition. 
\begin{definition}\label{S5.D1}
	Let $X,\, \partial X,\, U$ be Banach spaces and let the assumptions {\bf (H1)}--{\bf (H2)} be satisfied. The input-output system (\ref{input-output}) is called well-posed if for some $ \t>0 $ (hence all) there exists $ c_{\t}>0 $ such that the following inequality holds for all solutions of (\ref{input-output}):
	\begin{equation}\label{S5.2}
		\Vert z(\t)\Vert^{p}_{X}+\Vert y(\cdot)\Vert^{p}_{ L^{p}([0,\t];\partial X)}\leq c_{\t}\left( \Vert x\Vert^{p}_{X}+\Vert v(\cdot)\Vert^{p}_{ L^{p}([0,\t];\partial X)}+\Vert u(\cdot)\Vert^{p}_{ L^{p}([0,\t];U)}\right).
	\end{equation}
\end{definition}

In view of the above definition, we then obtain the following simple characterization of well-posed boundary input-output positive systems.
\begin{proposition}\label{S5.wellposed}
	Let $X, \partial X, U$ be Banach lattices and let the assumption {\bf ($\mathbf{H1}$)}-{\bf ($\mathbf{H2}$)} be satisfied. Furthermore, assume that $\mathbb{T},\,\Gamma,\,K$ are positive and $D_{\mu}$ is positive for every $\mu >s(A)$. Then, the system (\ref{input-output}) is well-posed if:
	\begin{itemize}
		\item[(\emph{i})] $ B$ is a positive $L^p$-admissible control operator for $A$, i.e., for some (hence all) $ \t>0 $, $$\Phi_{\tau}^Av:=\int_{0}^{\t}T_{-1}(\t-s)Bv(s)ds\,\in X_+,$$
		for all $v\in L^{p}_+(\R_+;\partial X) $.
		\item[(\emph{ii})] $C:=\Gamma_{\vert D(A)}$ is a positive $L^p$-admissible observation operator for $ A $, i.e., for some (hence all) $\alpha> 0$, $$\int_0^\alpha \Vert CT(t)x\Vert^{p} dt\leq  \gamma^{p}\Vert x\Vert^{p},$$  for all $ 0\leq x\in D(A)$ and a constant $\gamma:= \gamma(\alpha)>0$.
		\item[(\emph{iii})] For $ \tau >0 $, there exits a constant $ \kappa:=\kappa(\tau)>0 $ such that 
		\begin{equation}\label{wellposedestimate}
			\Vert \mathbb{F} v\Vert_{L^{p}([0,\tau];\partial X)}\leq \kappa \Vert v\Vert_{L^{p}([0,\tau];\partial X)},
		\end{equation}
		for $0 \le v\in W^{1,p}_{0}([0,\t];\partial X):=\left\{ v\in W^{1,p}([0,\t];\partial X ):v(0)=0\right\}$, where 
		\begin{equation}\label{wellposedABC}
			(\mathbb{F}v)(t):=\Gamma \Phi_{t}v,\qquad 0 \le v\in W^{1,p}_{0}([0,\t],\partial X),\;{\rm a.e.}\;t\in [0,\t].
		\end{equation}
	\end{itemize}
\end{proposition}
To proof the above proposition, we need the following lemma.
\begin{lemma}\label{S4.L2}
	Let $X,\partial X$ be Banach lattices, let $ A $ be a densely defined resolvent positive operator and $B\in \calL(\partial X,X_{-1})$. Define the vector space $Z\subset X$ by 
	$$
	Z:=R(\mu,A)\left(X+B \,\partial X\right),\;\; \mu> s(A). 
	$$
	If $B\in \calL(U,X_{-1})$ is a positive control operator, then $Z_1$ endowed with the norm
	$$
	\Vert z\Vert_{Z}^2=\inf\big\{\Vert x\Vert^2+\Vert v\Vert^2: x\in X,\,v\in \partial X,\;z=R(\mu,A)(x+Bv) \big\},
	$$
	is a Banach Lattice. 
\end{lemma}
{\it Proof}
It is clear that the definition of $Z$ is independent of the choice of $\mu> s(A)$. Moreover,  according to Lemma 4.3.12 (ii) of \cite{Sta}, $(Z,\Vert \cdot\Vert_{Z})$ is a Banach space satisfying $$ Z \hookrightarrow X.$$  
Then, it remains to prove that $Z$ is a vector lattice and $\Vert \cdot\Vert_{Z}$ is a lattice norm. Indeed, let $\mu> s(A) $, let $x\in X$ and $z\in Z$ such that $\vert x\vert \leq \vert z\vert $. Then $\vert x\vert  \leq R(\mu,A) \vert z_1\vert+R(\mu,A_{-1})B\vert v_1\vert $, since $z= R(\mu,A)z_1+ R(\mu,A_{-1})Bv_1$ for some $z_1\in X$ and $v_1\in \partial X$. By virtue of the decomposition property (see e.g \cite[Proposition II.1.6]{Schaf}), there exist $ x_1\in R(\mu,A)([-z_1,z_1])$ and  $ x_2\in R(\mu,A_{-1})B([-v_1,v_1])$ satisfying $x=x_1+x_2$. Therefore there exist $y_1\in [-z_1,z_1] $ and $y_2\in [-v_1,v_1] $ such that $x_1=R(\mu,A)y_1$ and $ x_2=R(\mu,A)By_2$. It follows that $x\in Z$ and hence $Z$ is a vector sublattice (as ideals are automatically lattice subspaces). Moreover, the fact that $B$ is positive and the norm on $X,\partial X$ are lattice norms yield that $\Vert \cdot\Vert_{Z}$ is a lattice norm. 
%	If $x\in E_+$ the then the symmetric order interval $[-x,x]$ is solid.
\qed

{\it Proof of Proposition \ref{S5.wellposed}}
Let $A:=(A_m)_{\vert \ker G}$ generate a positive C$_0$-semigroup $\mathbb{T}$ on $X$ and assume that  the operators $G,\Gamma,K$ are positive such that the operator $G$ is surjective. In view of Definition \ref{S5.D1}, we have to verify that the estimate (\ref{S5.2}) holds. According to \cite[Section 4]{El}, this is the same as characterizing the operators $A_m,G,\Gamma$ for which $(A,B,C)$ is a positive $L^p$-well-posed triplet on $X,\partial X, \partial X$. Since $B$ is a positive $L^p$-admissible control operator for $A$, then $\Phi_\t^A  L^{p}_+([0,\infty);\partial X)\subset X_+$ for all $\t\ge 0$. So, assuming that (without loss of generality) $0\in \rho(A)$ and using an integration by parts argument, one can show that
$$
\Phi_\t^A v =D_0 v(\t)-\Phi_\t\dot{v}\in Z_+,
$$
for all $0 \le v\in W^{1,p}_0([0,\t],\partial X)$ and $\t\geq 0$. In particular, we get that for every $\t\geq 0$ we have ${\rm Range}\,\Phi_\t^A\subset Z $, since $Z$ is a Banach lattice (Lemma \ref{S4.L2}). Therefore, the operator $\mathbb{F}$ is well-defined. Now, using a reasoning analog to \cite[Proposition 2.2]{El}, we show that the estimate (\ref{wellposedestimate}) yields that $(A,B,C)$ is a positive $L^p$-well-posed triplet on $X,\partial X, \partial X$. This ends the proof.
\qed

\begin{remark}\label{Regularity}
	We underline that the extended operator $\mathbb{F}$ is positive and bounded on $\calL(L^p([0,\t];\partial X))$ for each $\t\ge 0$. In particular, the extended output function $y$ of the system (\ref{input-output}) satisfies
	$$
	0 \le y(t;x,v)=C_{\Lambda} T(t)x+(\mathbb{F} v)(t),\qquad {\rm and\;a.e. }\; t\ge 0.
	$$
	for all $(x,v)\in X_+\times L_{loc,+}^{p}([0,\infty);\partial X)$, where $C_{\Lambda}$ is the \emph{Yosida extension} of $ C $ with respect to $ A $ and its domain denoted by $ D(C_{\Lambda}) $, consists of all $ x\in X $ for which $\underset{\mu \mapsto \infty}{\lim}C\mu  R(\mu ,A)x$ exists. The operator $\mathbb{F}$ is called the extended input-output control operator of $(A,B,C) $. Moreover, there exist $ \alpha> s(A) $ and a unique bounded analytic function $ \mathbf{H}:\mathbb{C}_{\alpha}\to \calL(U,Y) $ such that 
	\begin{equation}\label{transferfunction}
		%	\begin{array}{lll}
			\hat{y}(\mu)=CR(\mu,A)x+\mathbf{H}(\mu) \hat{v}(\mu),
			%	\end{array}
	\end{equation}
	for any $(x,v)\in X\times L^{p}_{loc}(\R_+;\partial X)\cap L^{p}_{\mu}(\R_+;\partial X)$ and $\mu\in \mathbb{C}_{\alpha}:=\{\mu\in \mathbb{C}:\;{\rm Re}\,\mu> \alpha\} $. Here $ \mathbf{H} $ denotes the transfer function of $(A,B,C) $ (or $\mathbb{F}$), see \cite[Chap. 4]{Sta} for more details. If, in addition, we have
	\begin{itemize}
		\item[\emph{(i)}] $\lim_{{\rm Re}\mu\to +\infty } y^*\mathbf{H}(\mu)v=0$ for all $v\in \partial X $ and $y^*\in (\partial X )'$, then $(A,B,C) $ is a positive $L^p$-well-posed weakly regular triplet (with feedthrough zero), where $(\partial X )^*$ is the dual of $\partial X $;
		\item[\emph{(ii)}] $\lim_{{\rm Re}\mu\to +\infty }\mathbf{H}(\mu)v=0$ in $\partial X $ for all $v\in \partial X $, then $(A,B,C) $ is a positive $L^p$-well-posed strongly regular triplet (with feedthrough zero).
	\end{itemize}
	According to \cite[Theorem 4.9]{El}, we have strong regularity and weak regularity of a positive $L^p$-well-posed triplet are equivalent. In particular, we simply say that $(A,B,C)$ is a positive regular triplet. In this case, we have $0 \le \Phi_t v\in D(C_\Lambda)$ and $0 \le (\F v)(t)=C_\Lambda \Phi_t v$ for any $v\in L^p_{loc,+}(\R_+,\partial X )$ and a.e. $t\ge 0$. In particular, we have that the state trajectory and the output function of (\ref{input-output}) satisfy $0\le z(t;x,v)\in D(C_\Lambda)$ and $0 \le y(t;x,v)=C_\Lambda z(t;x,v)$ for any $x\in X_+,$ $v\in L^p_{loc,+}(\R_+,\partial X )$ and a.e. $t\ge 0$, see \cite[Section 4]{El}.
\end{remark}
We also recall the following concept of positive admissible feedback operator, see \cite[Lemma 4.1]{El}.
\begin{definition}
	Let $X,\partial X$ be Banach lattices, let $(A,B,C)$ be a positive $ L^p $-well-posed regular triplet on $X,\partial X,\partial X$. We say that the identity $I_{\partial X}$ is a positive admissible feedback operator for $ (A,B,C) $ if and only if $r(\mathbb{F})<1 $.
\end{definition}

We end this section by the well-posedness of (\ref{S5.Sy1}) when the operators involved in (\ref{S5.Sy1}) are positive.
\begin{theorem}\label{wellposedness}
	Let $X,\, \partial X,\, U$ be Banach lattices and let the assumptions {\bf ($\mathbf{H1}$)}, {\bf ($\mathbf{H2}$)} and {\bf ($\mathbf{H3}'$)} be satisfied.	
	Then the operator $\mathcal{A}:D(\mathcal{A})\to X$ defined by
	$$
	\mathcal{A}=A_{m}, \qquad D(\mathcal{A})=\{x\in D(A_m):\: G x=\Gamma x\},
	$$
	generates a positive C$_0$-semigroup $\calT$ on $ X $. Moreover, if $K\in  \calL_+(U,\partial X)$, then the system (\ref{S5.Sy1}) has a unique mild solution $z$ satisfying
	\begin{equation}\label{variation}
		\begin{array}{lll}
			0\leq z(t)&=\calT(t)x+\displaystyle\int_{0}^{t}\calT_{-1}(t-s)BKu(s)ds\cr
			&:=\calT(t)x+\Phi_t^{\mathcal{A}}Ku,
		\end{array}
	\end{equation}
	for all $t\geq 0$, $ x\in X_+$ and $u\in L^p_+(\R_+;U)$. In particular, the perturbed boundary value control system \eqref{S5.Sy1} is positively well-posed. Furthermore, for $\mu >s(\calA)$,  
	$$
	\hat{z}(\mu)=R(\mu,\mathcal{A})x+\widehat{(\Phi_\cdot^\mathcal{A} u)}(\mu),\quad {\rm with }\quad \widehat{(\Phi_\cdot^\mathcal{A} u)}(\mu)=R(\mu,\mathcal{A}_{-1})B K\hat{u}(\mu), 
	$$
	for all $u\in L^p_\mu(\R_+;U) $, where $\hat{u}$ denote the Laplace transform of $u$.
\end{theorem}
{\it Proof}  By using the representation (\ref{S5.5}), the boundary input-output system (\ref{input-output}) can be reformulated as the following distributed-parameter system
\begin{equation}\label{ABC}
	\left\lbrace
	\begin{array}{lll}
		\dot{z}(t) =A_{-1} z(t)+Bv(t)+BKu(t),\qquad t> 0,\; z(0)=x,\\
		y(t) = C z(t), \qquad t> 0,
	\end{array}
	\right.
\end{equation}
where we recall that $C:=\Gamma_{\vert_{D(A)}}$ and $ B:=(\mu -A_{-1})D_\mu$ $( \mu > s(A)$). Therefore, according to \cite[Theorem 4.2]{El}, the system (\ref{ABC}) with the feedback law \eqref{feedback} is equivalent to the following open-loop system
$$	\left\lbrace
\begin{array}{lll}
	\dot{z}(t) =(A_{-1}+BC_\Lambda)z(t)+BKu(t),\qquad t> 0,\\
	z(0) =x, \qquad
\end{array}\right.
$$
since $ (A,B,C) $ is a positive $L^p$-well-posed regular triplet with the identity $ I_{\partial X} $ as a positive admissible feedback operator. In view of \cite[Theorem 4.2]{El}, the perturbed boundary control system (\ref{S5.Sy1}) has a unique mild solution $z$ satisfying 
\begin{equation}
	\begin{array}{lll}
		0\leq z(t)\in D(C_\Lambda), \qquad \forall\, x\in X_+,\; {\rm and \; a.e}\;\;t\ge 0;\\
		0\leq z(t)=\calT(t)x+\int_{0}^{t}\calT_{-1}(t-s)BKu(s)ds,
	\end{array}
\end{equation}
for all $t\geq 0$, $ x\in X_+$ and $u\in L^p_+(\R_+;U)$, where $\calT$ is the positive C$_0$-semigroup generated by $(A_{-1}+BC_\Lambda)$. Moreover, according to \cite[Theorem 4.3]{El}, we have $\mathcal{A}=(A_{-1}+BC_\Lambda)$. The last claim follows by applying Laplace transform to both side of (\ref{variation}). This completes the proof.
\qed

	\section{Proof of the main results }\label{Sec3}
In this section, we give the proof of the main results of this paper. Firstly, we recall that the system (\ref{S5.Sy1}) is given by 
$$
\left\lbrace
\begin{array}{lll}
	\dot{z}(t)= A_{m}z(t),\qquad \quad\qquad t> 0,\\
	z(0)=x,\\
	(G-\Gamma)z(t)=Ku(t), \qquad t> 0.
\end{array}
\right.
$$
Moreover, to (\ref{S5.Sy1}) we associate the following linear (differential) operator 
\begin{equation*}
	\mathcal{A}:= A_m, \qquad D(\mathcal{A})= \left\{x\in D(A_m):\quad Gx= \Gamma x\right\}.
\end{equation*}
\subsection{Proof of Theorem \ref{Main-T1}}
Under the assumptions {\bf(A1)}-{\bf(A3)} and according to \cite[Theorem 4.3]{HMR}, the operator $\calA$ generate a C$_0$-semigroup $\calT$ on $X$ and the state trajectory of (\ref{S5.Sy1}) satisfies the following variation of constant formula
$$
z(t;x,u)=\calT(t)x+\Phi_t^\mathcal{A}K u,\qquad t\geq 0,
$$
for all $x\in X$ and $u\in L_{loc}^p(\R_+;U)$, where we recall that 
$$
\Phi_t^\mathcal{A}K u=\int_{0}^{t}\calT_{-1}(t-s)BKu(s)ds.
$$
So, given a prescribed time $\t>0$, we shall be concerned with the final state
$$
z(\t;x,u)=\calT(\t)x+  \Phi_\t^\mathcal{A} Ku, 
$$
and the following space of reachable states from the origin with respect to positive controls in time $\t$
$$
\calR_\t^{U_+}:= \left\{z(\t;0,u)\,: u\in L_+^p([0,\t];U)  \right\}.
$$
Then, the concept of boundary approximate controllability with respect to positive controls in finite time is defined as follows.
\begin{definition}\label{definition1}
	Let the assumptions of Theorem \ref{Main-T1} be satisfied. We say that the system (\ref{S5.Sy1}) is boundary approximately controllable with respect to positive controls if the reachable set from the origin in finite time %$\mathcal{R}^{U_+}$  defined by 
	\begin{align*}
		\mathcal{R}^{U_+}&:=\bigcup_{\t>0}\mathcal{R}_\t^{U_+}\\
		&=\left\{z\in X\,\vert \exists\,\t_0>0\; {\rm and}\; u\in L^p_+([0,\t_0];U)\;{\rm such\; that }\; z= \Phi_{\t_0}^\mathcal{A} Ku\right\},
	\end{align*}
	is dense in $X$.
\end{definition}

{\it Proof Theorem \ref{Main-T1}}
In view of Definition \ref{definition1} the system (\ref{S5.Sy1}) is boundary approximately controllable with respect to positive controls if and only if $\calR^{U_+}$ is dense in $X$. Notice that $\calR^{U_+}$ is a convex cone with $0\in \calR^{U_+} $. 

Let $\mu\ge\omega>\omega_{0}(\calA)$ and $\varphi\in X'$. Assume that $\overline{\calR^{U_+}}=X$ and
\begin{equation}\label{polar-posit}
	\langle (I_{ X}- D_{\mu}\Gamma)^{-1}D_{\mu}K \hat{u}(\mu),\varphi\rangle \leq 0, \qquad\qquad \forall\, u \in L^p_{\mu,+}(\R_+;U).
\end{equation}
By the uniqueness of the Laplace transform (\cite[Theorem 1.7.3]{ACMF}), one can see that (\ref{polar-posit}) implies that 
\begin{equation}\label{positive}
	\langle  \int_{0}^{\t}\calT_{-1}(\t-s)BKu(s)ds,\varphi\rangle \leq 0, \qquad \forall\, u\in L^p_{\mu,+}(\R_+;U),\,\t> 0,
\end{equation}
since 
$$
R(\mu,\mathcal{A}_{-1})z= (I_X-D_{\mu}\Gamma)^{-1} R(\mu, A_{-1})z, \qquad \forall\, z\in X_{-1}.
$$
Therefore, by virtue of Lemma \ref{Main-lemma}, we get that $\varphi =0$ (as $\overline{\calR^{U_+}}=X$). 

Conversely, let us assume that there exists $\omega\in \R$ such that $1\in \rho(\Gamma D_\mu)$ and (\ref{charactera-T1}) holds for all $\mu\ge \omega$. Furthermore, we assume that $ \overline{\calR^{U_+}}\neq X$. Then, according to Lemma \ref{Main-lemma}, there exists $0\neq \varphi\in X'$ such that (\ref{positive}) holds. Taking the Laplace transform in (\ref{positive}), we get
$$
\langle R(\mu,\mathcal{A}_{-1})B K \hat{u}(\mu),\varphi\rangle \leq 0 ,$$
or equivalently,
$$ \left\langle (I_{ X}- D_{\mu}\Gamma)^{-1} R(\mu, A_{-1}) B K \hat{u}(\mu),\varphi\right\rangle \leq 0.
$$
It follows that $\varphi\in \left((I_{X}-D_{\mu}\Gamma)^{-1}D_{\mu}KU_{+}\right)^{\circ} $-- contradiction--. Furthermore, for $\omega\in \R$ such that $1\in \rho(\Gamma D_\mu)$ for all $\mu \ge \omega$, we have 
$$
\begin{array}{lll}
	(I_X-D_{\mu}\Gamma)^{-1} R(\mu, A_{-1}) B Ku
	&=(I_X-D_{\mu}\Gamma)^{-1} D_\mu Ku\\
	&=D_\mu(I_{\partial X}-\Gamma D_{\mu})^{-1} Ku,
\end{array}
$$
for all $u\in U_+$ and $\mu \ge \omega$. This completes the proof.
\qed

\subsection{Proof of Theorem \ref{Main-T2}}
Here, we assume that the operators involved in (\ref{S5.Sy1}) are all positive. As such, we focus our attention on its boundary approximate controllability with respect to positive controls. As the states are confined to the positive orthant, one need to consider another concept of approximate controllability of (\ref{S5.Sy1}), namely approximate positive controllability introduced in \cite{Tilman}. Indeed, let us consider the following space of reachable positive states from the origin in time $\t$ with respect to positive controls
$$
\calR_\t^+:= \left\{z(\t;0,u)\,: u\in L^p_+([0,\t];U)  \right\}.
$$
Then the concept of boundary approximate positive controllability in finite time is defined as follows.
\begin{definition}\label{definition2}
	Let the assumptions of Theorem \ref{Main-T2} be satisfied. We say that the system (\ref{S5.Sy1}) is boundary approximately positive controllable if the reachable set from the origin in finite time
	\begin{align*}
		\mathcal{R}^{+}&:=\bigcup_{\t>0}\mathcal{R}_\t^{+}\\
		&=\left\{z\in X_+\,\vert \;\exists\,\t_0>0\; {\rm and}\; u\in L^p_+([0,\t_0];U)\;{\rm such\; that }\; z= \Phi_{\t_0}^\mathcal{A} Ku\right\},
	\end{align*}
	is dense in $X_+$.
\end{definition}

{\it Proof of Theorem \ref{Main-T2}}
It follows from Theorem \ref{wellposedness} that the system (\ref{S5.Sy1}) has a unique positive mild solution $z$ satisfying the following variation of constant formula
$$
0\le z(t)=\calT(t)x+\Phi_t^{\mathcal{A}}Ku,
$$
for all $t\geq 0$, $ x\in X_+$ and $u\in L^p_{+}(\R_+;U)$, where we recall that $\calT$ is the C$_0$-semigroup generated by $\calA$.

Now, let $\mu > \max\{\omega,s(\mathcal{A})\}$ and $ \varphi\in X'$. Assume that (\ref{S5.Sy1}) is boundary approximately positive controllable and $\langle (I_{\partial X}- D_{\mu}\Gamma)^{-1}D_{\mu}K \hat{u}(\mu),\varphi\rangle \leq 0$ for all $u \in L^p_{\mu,+}(\R_+;U)$. The uniqueness of the Laplace transform yields that
$$
\left\langle  \int_{0}^{\t}\calT_{-1}(\t-s)BKu(s)ds,\varphi\right\rangle \leq 0,
$$
for all $u \in L^p_{\mu,+}([0,\infty);U)$ and $\t>0$, since 
$$
R(\mu,\mathcal{A}_{-1})z= (I_X-D_{\mu}\Gamma)^{-1} R(\mu, A_{-1})z, \qquad \forall\, z\in X_{-1,+}\;(:=X_+\cap X_{-1}).
$$
It follows that $\varphi \le 0$ as $\overline{\calR^{+}}=X_+$).

Conversely, let us assume that there exits $\omega>s(A)$ such that $r(\Gamma D_{\omega})<1$ and the following implication holds for all $\varphi\in X'$:
$$
\langle (I_{ X}- D_{\mu}\Gamma)^{-1}D_{\mu}K u,\varphi\rangle \leq 0, \qquad \forall\, u \in U_+,\;\mu > \max\{\omega,s(\mathcal{A})\}\Longrightarrow \;\varphi\le 0.
$$
Moreover, assume that $\overline{\calR^{+}}\neq X_+$. Then, by Hahn-Banach Theorem there exist $z_1\in X_+ \setminus \calR^{+} $, $0\neq \varphi\in X'$ and $\beta\in \R$ such that 
\begin{equation}\label{Dada}
	\langle \Phi_\t^{\mathcal{A}}Ku,\varphi\rangle  < \beta < \langle z_{1} ,\varphi\rangle,
\end{equation}
for all $u \in L^p_{\mu,+}(\R_+;U)$ and $\t>0$. On the other hand, let $v\in U_+$ and define the following sequence of functions: 
$$
u_{n}(t):=
\left\lbrace
\begin{array}{lll}
	0	, &{\rm if}\;\; 0\leq t \leq  \frac{1}{n},
	\\
	nv, &{\rm if}\;\; \frac{1}{n}< t .
\end{array}
\right.
$$
Clearly, we have $u_n\in L^p_{\mu,+}(\R_+;U)$. So, for $\mu>\max\{\omega,s(\calA)\}$ and taking Laplace transform in (\ref{Dada}), 
$$
\langle  (I_{ X}- D_{\mu}\Gamma)^{-1}D_{\mu}K\hat{u}_{n}(\mu),\varphi\rangle < \beta < \langle z_{1} ,\varphi\rangle.
$$
Thus,
$$
\left\langle  (I_{ X}- D_{\mu}\Gamma)^{-1}D_{\mu}Kn\int_{\frac{1}{n}}^{\infty}e^{-\mu t}vdt,\varphi\right\rangle < \beta
$$
and hence,
$$
\left\langle  (I_{ X}- D_{\mu}\Gamma)^{-1}D_{\mu}K\int_{\frac{1}{n}}^{\infty}e^{-\mu t}vdt,\varphi \right\rangle < \frac{1}{n}\beta ,
$$
with the right-hand side converging to $0$ as $n\to \infty$. Then,
$$
\langle  (I_{ X}- D_{\mu}\Gamma)^{-1}D_{\mu}K v,\varphi\rangle < 0.
$$
Therefore, $ \varphi\in\left((I_{ X}- D_{\mu}\Gamma)^{-1}D_{\mu}KU_{+}\right)^{\circ}  $ and $0 < \langle z_{1} ,\varphi\rangle$ (as $0\in \calR^+$), a contradiction.  Finally, for $\omega>s(A)$ such that $r(\Gamma D_\omega)<1$, we have 
$$
\begin{array}{lll}
	(I_X-D_{\mu}\Gamma)^{-1} R(\mu, A_{-1}) B Ku
	&=D_\mu(I_{\partial X}-\Gamma D_{\mu})^{-1} Ku\\
	&=D_{\mu}\sum_{n=0}^{\infty} (D_{\mu} \Gamma )^{n}Ku,
\end{array}
$$
for all $u\in U_+$ and $\mu\ge \omega$, where we have used the Neumann series representation of $(I_{\partial X}-\Gamma D_{\mu})^{-1}$. Thus,
$$
\left((I_{ X}- D_{\mu}\Gamma)^{-1}D_{\mu}KU_{+}\right)^{\circ}=\bigcap\left\{\big(D_{\mu}(\Gamma D_{\mu})^n KU_{+}\big)^{\circ},\;\; n\in \N \right\},
$$
and this ends the proof.
\qed

	\section{Application}
In this section we illustrate our abstract results through two relevant applications of PDEs on networks. Such evolutionary systems on networks have been studied by many authors, in particular we refer to the monograph \cite{Delio} and the proceedings \cite{Ali}.

\subsection{Flows on closed networks}
Let us consider the following system of transport equations on a network:
$$
(\Sigma_{\mathsf{TN}})	
\left\lbrace	
\begin{array}{llll}
	\frac{\partial}{\partial t} z_{j}(t,x,v)= v\frac{\partial }{\partial x}z_{j}(t,x,v)+q_{j}(x,v).z_{j}(t,x,v),\quad t\geq 0, \,(x,v)\in \Omega_j,\\  
	z_{j}(0,x,v)= f_{j}(x,v)\ge 0, \qquad \qquad \qquad \qquad \qquad \quad \; (x,v)\in \Omega_j, \;({\rm IC}) \\
	\imath^{out}_{ij}z_{j}(t,1,\cdot)= \mathsf{w}_{ij}\displaystyle\sum_{k\in \calM} \imath^{in}_{ik}\mathbb{J}_k(z_{k})(t,0,\cdot)+\sum_{l\in \calN_c}\mathsf{b}_{il}u_{l}(t,.), \; t\geq 0,\;({\rm BC})
\end{array}
\right.
$$
for $ i\in \{1,\ldots,N\}:=\calN$, $j\in \{1,\ldots,M\}:=\calM$ and $l\in \{1,\ldots,n\}:=\calN_c$ with $M\ge N\ge n$, where we set $\Omega_j:=[0,1]\times [v_{\min},v_{\max}]$. The corresponding transport equations are defined on the edges of a finite graph $\mathsf{G}$, whose edges are normalized and identified with the interval $[0,1]$ with endpoints "glued" to the graph structure. The connection of such edges being described by the coefficients $\imath^{out}_{ij},\imath^{in}_{ik}\in \{0,1\}$. The flow velocity along the edges is determined by the function $v$, whereas its absorption is determined by the functions $ q_j(\cdot, \cdot) $. The boundary condition (BC) determines the propagation of the flow along the various components of the network. The weights $0\le \textsf{w}_{ij}\le 1$  express the proportion of mass being redistributed into the edges and the non-local operators $\mathbb{J}_k$ describe the scattering at the vertices. Moreover, for $i,l\in  \calN\times \calN_c$, the coefficients $\mathsf{b}_{il}\ge 0$ denotes the entries of the input matrix $ K $ and $ u_l\ge 0 $ denotes the input functions at the vertices. The system $(\Sigma_{\mathsf{TN}})$ is a generalization of the transport network in (\cite{EHR1}).

Next, we are concerned with the well-posedness and boundary approximate controllability with respect to positive controls of $(\Sigma_{\mathsf{TN}})$. To this end, we need to recall some notation from graph theory. Here and in the following, we consider a finite connected metric graph $\mathsf{G}=(\mathsf{V},\mathsf{E})$ and a flow on it (the latter is described by $(\Sigma_{\mathsf{TN}})$). The graph $ \mathsf{G} $ is composed by $ N\in \mathbb{N} $ vertices $ \alpha_{1},\,\ldots,\alpha_N $, and by $ m\in \mathbb{N} $ edges $ \mathsf{e}_{1},\,\ldots,\mathsf{e}_M $ which are normalized so as to be identified with the interval $ [0,1] $. Each edge is parameterized contrary to the direction of the flow on them, i.e., the material flows from $1$ to $0$. The topology of the graph $\mathsf{G}$ is described by the incidence matrix $\mathcal{I}=\mathcal{I}^{out}-\mathcal{I}^{in}$, where $\mathcal{I}^{out}$ and $\mathcal{I}^{in}$ are \emph{the outgoing incidence matrix} and the \emph{incoming incidence matrix} of $\mathsf{G}$ having entries
$$
\imath^{out}_{ij}:=
\left\lbrace	
\begin{array}{lll}
	1,\quad {\rm if  } \begin{tikzpicture}
		\node (P) at (0,0) {$\mathsf{v}_i$};
		\node (T) at (0.7,0.2) {$\mathsf{e}_j$};
		\node (S) at (1.2,0) {};
		\draw[*->,>=latex] (P) to[=1] (S);
	\end{tikzpicture},
	\\
	0,\quad {\rm if\, not},
\end{array}\right.\quad
\imath^{in}_{ij}:=		\left\lbrace	
\begin{array}{lll}
	1, \quad {\rm if  }  \begin{tikzpicture}
		\node (P) at (0,0) {};
		\node (T) at (0.4,0.2) {$\mathsf{e}_j$};
		\node (S) at (1.2,0) {$\mathsf{v}_i$};
		\draw[->*,>=latex] (P) to[=6] (S);
	\end{tikzpicture},
	\\
	0, \quad {\rm if\, not},
\end{array}\right.
$$
respectively. Replacing $1$ by $\textsf{w}_{ij}\geq 0$ in the definition of $\imath^{out}_{ij}$, we obtain the so-called weighted outgoing incidence matrix $\mathcal{I}^{out}_{\mathsf{w}}:=(\mathsf{w}_{ij}\imath^{out}_{ij})$. In this cases, $\mathsf{G}$ is called a weighted graph and its topology is described via the weighted transposed adjacency matrix $\mathbb{A}:=\mathcal{I}^{in}(\mathcal{I}_{\mathsf{w}}^{out})^{\top}$ given by, for $i,k\in \calN$,
\begin{equation}\label{weighted}
	(\mathbb{A})_{ik}:= 		\left\lbrace	
	\begin{array}{lll}
		\mathsf{w}_{kj},& {\rm if}\; \exists\, \mathsf{e}_j\begin{tikzpicture}
			\node (P) at (2,0) {};
			\node (R) at (2.26,0.2) {$\mathsf{v}_i$};
			\node (T) at (1.3,0.2) {$\mathsf{e}_j$};
			\node (S) at (4,0) {};
			\node (P') at (0.5,0) {};
			\node (P'') at (2.26,0) {};
			\node (T') at (0.5,0.2) {$\mathsf{v}_k$};
			\draw[*->*,>=latex] (P') to[=0] (P'');
		\end{tikzpicture},
		\\
		0, & {\rm if\, not}.
	\end{array}\right.
\end{equation}

In the rest of this section, for $ 1\le p<+\infty $, let us consider the Banach spaces $(Y_p,\|\cdot\|_{Y_p})$, $(X_p,\|\cdot\|_{X_p})$ and $(\mathbb{W}_p,\Vert \cdot \Vert_{\mathbb{W}_p})$ defined by
$$
\begin{array}{llll}
	Y_p&:=L^{p}([v_{\min},v_{\max}])^M,\qquad \Vert f \Vert^{p}_{Y_p} :=\sum_{j=1}^{M}\Vert f_j\Vert^{p}_{L^{p}([v_{\min},v_{\max}])},\\
	X_p&:=L^{p}([0,1],Y_p), \qquad\qquad \|\varphi \|^{p}_{_p} :=\int_{0}^{1} \|\varphi(x,\cdot)\|_{Y_p}^{p} dx,\\
	\partial X_p&:= L^p([v_{\min},v_{\max}])^N \qquad \Vert g \Vert^{p}_{\partial X_p} :=\sum_{j=1}^{N}\Vert f_j\Vert^{p}_{L^{p}([v_{\min},v_{\max}])},\\
	\mathbb{W}_p&:= W^{1,p}([0,1];Y_p),\qquad \Vert f \Vert_{\mathbb{W}_p}^p=\Vert f\Vert^{p}_{X_p}+ \Vert \partial_x f\Vert^{p}_{X_p}.
\end{array}
$$
Moreover, let us consider the following operators 
\begin{equation}\label{S1.Am}
	\begin{array}{llll}
		A_{m}f&= v\partial_x f+q(\cdot,\cdot)f,\\
		D(A_{m})&= \left\{ f =(f_{j})_{j\in \calM}\in \mathbb{W}_p:  f(1,v)\in {\rm Range }\, (\mathcal{I}^{out}_{\mathsf{w}})^{\top}\right\},
	\end{array}
\end{equation}
where $q(\cdot,\cdot):={\rm diag}\,(q(\cdot,\cdot))_{j\in \calM}$.  Moreover, we set $\mathbb{J}={\rm diag}(\mathbb{J}_k)_{k\in \calM}$ where the scattering operators $\mathbb{J}_k$ are given by 
\begin{equation}
	\mathbb{J}_k(f_k)(x,v)=\int_{v_{\min}}^{v_{\max}}\ell_k(x,v,v')f_k(x,v') dv',\qquad (x,v)\in \Omega_k,\; f\in X_p,
\end{equation} 
where $ \ell_k\in L^{\infty}_+(\Omega_k\times [v_{\min},v_{\max}])$ for every $k\in \calM$. We select 
\begin{equation}\label{Operators-boundary}
	G f:=	\mathcal{I}^{out} f(1,\cdot),\quad
	\Gamma f :=\mathcal{I}^{in}(\mathbb{J}f)(0,\cdot),\quad v\in [v_{\min},v_{\max}],\; f\in \mathbb{W}_p.
\end{equation}
We also introduce the input space $U:=L^p([v_{\min},v_{\max}])^n$ and the control operator $K $ is given by 
$$
(Ku)_i(t,\cdot) =\sum_{l\in \calN_c}\mathsf{b}_{il}u_l(t,\cdot) ,\qquad \forall\, i\in \calN,\,t\ge 0,
$$
where $\mathsf{b}_{il}\ge 0$ for all $i,l\in \calN\times \calN_c$.

With the above notation, one can rewrite the system $(\Sigma_{\mathsf{TN}})	$ as 
$$
\left\lbrace
\begin{array}{lll}
	\dot{z}(t)= A_{m}z(t),& t> 0,\\
	z(0)=x,\\
	(G-\Gamma)z(t)=Ku(t), & t> 0,
\end{array}
\right.
$$
where $z(t)=(z_j(t,\cdot,\cdot))_{j\in \calM},\,f=(f_j(\cdot,\cdot))_{j\in \calM}$ and $u(t)=(u_l(t,\cdot))_{l\in \calN_c}$.

In order to apply the results of the previous sections, let us introduce the following assumptions:
\begin{itemize}
	\item[{\bf (A1)}] $0<v_{\min}\le v\le v_{\max}$ and $ q_j(\cdot, \cdot)\in L^{\infty}(\Omega_j) $ for every $j\in \calM$.
	\item[{\bf (A2)}] Each vertex has at least one outgoing edge.
	\item[{\bf (A3)}] The weights $\textsf{w}_{ij}$ satisfy 	$\sum_{j\in \calM}\textsf{w}_{ij}=1,\; \forall\, i\in \calN$.
\end{itemize}
In the above, {\bf (A1)} specifies the transport process along each edges, {\bf (A2)} is equivalent to the statement that $\mathcal{I}^{out}$ is surjective, while {\bf (A3)} implies that the free delay boundary condition (BC) exhibits standard \emph{Kirchhoff conditions} and the matrix $\mathbb{A}$ is column stochastic. Note also that, under the condition {\bf (A1)} it is not difficult to prove that the operator $(A, \ker G)$ generates a strongly continuous positive semigroup $(T(t))_{t\geqslant 0} $ on $X_p$ given by
\begin{equation}\label{Semigroup}
	(T(t)f)_j(x,v)=	\left\lbrace
	\begin{array}{lll}
		e^{\int_0^{t}q_j(x+v\sigma,v)d\sigma}f_j(x+vt,v), & {\rm if }\;  x+vt \leq 1,
		\\
		0,& { \rm if\, not }.
	\end{array}
	\right.
\end{equation}
for all $f\in X_p$, $(x,v)\in \Omega_j$ and $j\in\calM$. 
\begin{lemma}\label{Transport-wellposed}
	Let the assumptions {\bf (A1)}-{\bf (A3)} be satisfied. Then the operator 
	\begin{equation}\label{CP1}
		\calA=A_{m},\quad
		D(\calA)=\big\lbrace f\in D(A_m):\; \mathcal{I}^{out}(\mathbb{J}f)(1,v)=\mathcal{I}^{in}(\mathbb{J}f)(0,v) \;\big\rbrace.
	\end{equation}
	generates a positive C$_0$-semigroup $\calT$ on $ X_p $. Moreover, for $\mu_0>s(A)$ such that $r(\mathbb{A}(\mu_0))<1$, we have
	\begin{equation}\label{3.32}
		R(\mu,\calA)=(I_{X_p}+D_\mu(I_{\partial X_p}-\mathbb{A}(\mu))^{-1}\Gamma)R(\mu,A),
	\end{equation}
	for all $\mu\ge \mu_0$, where we set $ \mathbb{A}(\mu):=\Gamma D_{\mu}$ with
	$$
	\begin{array}{lll}
		(D_\mu g)(x,v)&=e^{\int_{x}^{1}\frac{q_j(\sigma,v)-\mu}{v}d\sigma}\displaystyle\sum_{i\in \calN}\mathsf{w}_{ij}g_i(v),\\
		(R(\mu,A)f)_j(x,v)&=\int_{x}^{1}	e^{\int_0^{y}\frac{q_j(\sigma,v)-\mu}{v} d\sigma}vf_j(y,v)dy,
	\end{array}
	$$
	for all $g\in \partial X_p$, $(x,v)\in \Omega_j$ and $j\in \calM$.
\end{lemma}	
{\it Proof}
First, observe that $\mathbb{J}\ge 0$ as $ \ell_j\in L^{\infty}_+(\Omega_j\times [v_{\min},v_{\max}])$ for all $j\in \calM$. To prove the generation of $(\calA,D(\calA))$ on $X_p$ we shall use \cite[Theorem 4.3]{El}. In fact, let $A$ be the restriction operator of $A_m$ on $\ker G$, i.e. $A=(A_m)_{\vert \ker G}$. Then, by (\ref{Semigroup}), the operator $A$ is clearly densely defined resolvent positive. On the other hand, in view of the assumption {\bf (A1)}, the operator $G$ is surjective. Moreover, by a simple computation we get that the Dirichlet operator $D_\mu$ associated to $G$ is given by
$$
D_\mu={\rm diag}(e^{\int_{\cdot}^{1}\frac{q_j(\sigma,v)-\mu}{v}d\sigma})_{j\in \calM}\otimes(\mathcal{I}^{out}_{\mathsf{w}})^{\top}, \quad {\rm Re}\,\mu > \tilde{q}:=-\inf\Vert q_j\Vert_\infty,
$$
where $\otimes$ denotes the tensor product. Now, let us define the operators  $C:=\Gamma_{\vert \ker G}$ and $B=(\mu-A_m)D_\mu$ for $\mu >\tilde{q}$. Using the injectivity of Laplace transform one can show 
$$
\begin{array}{lll}
	(\Phi_t^Ag)_j(x,v)&:=\left(\int_{0}^{t}T_{-1}(t-s)Bg(s)ds\right)(x,v)\\
	&=\sum_{i\in \calN}	e^{\int_x^{1}\frac{q_j(\sigma,v)}{v}d\sigma}\mathsf{w}_{ji}g_i(\frac{tv-1+x}{v}) \1_{[\frac{1-x}{v},\infty)}(t)
\end{array}		
$$
for all $g\in L^p(\R_+;\partial X_p),\,t\ge 0,\,(x,v)\in \Omega_j$ and $j\in \calM$. So, for every $t\ge 0$ we have $\Phi_t^A\in \calL_+(L^p(\R_+;\partial X_p) ,X_p)$ and hence $B$ is a positive $L^p$-admissible control operator for $A$. Moreover, for $0<\alpha<\frac{1}{v_{\max}}$ we have 
$$
\int_{0}^{\alpha}\Vert CT(t)f\Vert_{\partial X_p}^p dt\le c_p N(v_{\max}-v_{\min})^p\sup_{j\in \calM}\Vert \ell_j\Vert^p_{\infty} e^{\left(p\sup_{j\in \calM}\Vert q_j\Vert_{\infty}\right)} \Vert f\Vert_{X_p}^p,
$$
for all $0\le f\in D(A)$, where we have used H\"{o}lder's inequality for $ \frac{1}{p}+\frac{1}{q}=1 $. Then, $C$ is a positive $L^p$-admissible observation operator for $A$. Furthermore, for $g\in W^{1,p}(\R_+,\partial X_p)$ with $g(0)=0=\dot{g}(0)$, we have
$$
\begin{array}{lll}
	(\mathbb{F}g)_k(t)&=(C\Phi_t^{A}g)_k\cr
	&= \displaystyle\sum_{j=1}^{M}\sum_{i=1}^{N}	\imath_{kj}^{in}\int_{v_{\min}}^{v_{\max}}\ell_j(0,v,v')e^{\int_0^{1}\frac{q_j(\sigma,v')}{v'}d\sigma}\mathsf{w}_{ji}g_i(\frac{tv'-1}{v'}) \1_{[\frac{1}{v'},\infty)}(t)dv',
\end{array}
$$
for $t\ge 0$ and $k\in \calN$. With this explicit expression of the input-output control operator $\mathbb{F}$ and according to \cite[Proposition 4.3]{El}, it is not difficult to see that $(A,B,C)$ is a positive $L^p$-well-posed triplet on $X_p,\partial X_p,\partial X_p $. In addition, for $\lambda>0$ and $\mu> \tilde{q}$ such that $\mu\neq \lambda$, using a computation involving the resolvent identity (we omit the details) we obtain 
$$
\lim_{\lambda\mapsto +\infty} \lambda CR(\lambda,A)D_\mu g=\Gamma D_\mu g, \qquad \forall\, g\in \partial X_p.
$$
It follows that $D_\mu \subset D(C_{\Lambda})$ and $(C_{\Lambda})_{\vert D(A_m)}=\Gamma$,  where $C_{\Lambda} $ is the \emph{Yosida extension} of $ C:=\Gamma_{\vert D(A)} $ with respect to $ A$. In particular, $(A,B,C)$ is a positive $L^p$-well-posed regular triplet. On the other hand, for $t< \tfrac{1}{v_{\max}}$ we have $I_{\partial X_p}-\mathbb{F}=I_{\partial X_p}$. Thus, according to \cite[Lemma 4.1.]{El}, $I_{\partial X_p}:\partial X_p\to \partial X_p$ is a positive admissible feedback operator for $(A,B,C)$. Therefore, according to \cite[Theorem 4.3]{El}, $\calA$ generates a positive C$_0$-semigroup $\calT$ on $X_p$. 
\qed

Under the assumptions of Lemma \ref{Transport-wellposed} and according to Theorem \ref{wellposedness}, we get that the transport network system $(\Sigma_{\mathsf{TN}})$ is positively well-posed and hence it has a unique positive mild solution $z(\cdot):\R_+\to X_p$ satisfying  
$$
\begin{array}{lll}
	0\le z(t)\in D(C_{\Lambda})&,\qquad {\rm for \; a.e.} \qquad t\ge 0,\cr
	0 \le z(t;f,u)=&\calT(t)f+\int^t_0 \calT_{-1}(t)(t-s)BKu(s)\,ds:=\calT(t)f+\Phi_{t}^{\calA}u,
\end{array}
$$
for all $t\ge 0$, $f\in (X_p)_+$ and $u\in L^p_+(\R_+,U)$. Moreover, for $\omega>0$ (large enough), we have $r(\mathbb{A}(\omega))< 1$. Then, according to Theorem \ref{Main-T2}, the system $(\Sigma_{\mathsf{TN}})$ is boundary approximately positive controllable if and only if,  
\begin{equation}\label{F2'}
	\bigcap\left\{\left(D_\mu\mathbb{A}^m(\mu) KU_{+}\right)^{\circ},\;\; m\in \N \right\} =(X_p)^{\circ}_+,  \qquad \forall\, \mu\ge \omega.
\end{equation}
Furthermore, if $q_j\equiv 0$, we obtain a more compact criteria as follows.
\begin{theorem}\label{Theorem-transport}
	Assume that $q_j\equiv 0$ for all $j\in \calM$ and let the assumptions {\bf (A1)}-{\bf (A3)} be satisfied. Then, the system $(\Sigma_{\mathsf{TN}})$ is boundary approximately positive controllable if and only if 
	\begin{equation}\label{F1}
		\bigcap\left\{\left( (\mathcal{I}^{out}_{\mathsf{w}})^{\top}\mathbb{A}^m (\mu) KU_{+}\right)^{\circ},\;\; m\in \N \right\} =(Y_p)^{\circ}_+,\quad \forall\, \mu\ge \mu_0,
	\end{equation}
	for some $\mu_0>0$ sufficiently large. Here, the operator $\mathbb{A}(\mu)$ is given by 
	$$(\mathbb{A}(\mu)g)_k(v)= \sum_{j\in \calM}\sum_{i\in \calN}	\imath_{kj}^{in}\int_{v_{\min}}^{v_{\max}}\ell_j(0,v,v')e^{-\frac{\mu}{v'}(1-x)}\mathsf{w}_{ji}g(v')\; dv',$$
	for all $k\in \calN$ and $g\in \partial X_p$.
\end{theorem}  %{\rm Re}\,
{\it Proof}
Note that, for $p,q\in [1,+\infty)$, we have 
$$
L^p([0,1],Y_p)=\overline{ L^p([0,1]) \otimes Y_p}\quad {and}\quad
L^q([0,1],Y_q)=\overline{L^q([0,1]) \otimes (Y_q)'},
$$
see e.g. \cite[Section 2.2.3]{Introduction-tensor-book}. Then, for $ \frac{1}{p}+\frac{1}{q}=1 $, we have 
\begin{equation}\label{vertor-space-equality}
	(X_p)^{\circ}_+:= -(X_p)_+'=\overline{L^q_+([0,1])\otimes  (Y_p)_+^{\circ} }.
\end{equation} 
To prove our claim we will use Theorem \ref{Main-T2}. Indeed, using the explicit expression of $\mathbb{A}(\mu)$ we get
$$
\Vert \mathbb{A}(\mu)\Vert \le \kappa e^{-\frac{\mu}{v_{\min}}},
$$
where $\kappa$ is a generic constant depending on $v_{\min},\;v_{\max}$, $\ell$ and the exponent $p$. Thus, $\Vert \mathbb{A}(\mu)\Vert \to 0$ as $\mu \mapsto\infty$. So, there exists $\mu_0>0$ large enough such that $\Vert \mathbb{A}(\mu)\Vert<1$ for all $\mu \ge \mu_0$. 

Now, let $\mu \ge \mu_0$ and $v\in [v_{\min}, v_{\max}]$ be fixed. It follows from (\ref{F2'}) and (\ref{vertor-space-equality}) that $(\Sigma_{\mathsf{TN}})$ is boundary approximately positive controllable if and only if
\begin{equation}\label{F3}
	\bigcap\left\{\left(D_\mu\mathbb{A}^m(\mu) KU_{+}\right)^{\circ},\;\; m\in \N \right\} =\left(L^q_+([0,1])\otimes  (Y_p)_+^{\circ} \right),  \qquad \forall\, \mu\ge \mu_0.
\end{equation}
The explicit expression of $D_\mu$ together with Lemma \ref{density1} further yield that $(\Sigma_{\mathsf{TN}})$ is boundary approximately positive controllable if and only if the condition (\ref{F1}) holds. This ends the proof. 
\qed

In particular, for a simple transport process we have the following Kalman-type rank condition.
\begin{corollary}\label{COR1}
	Let the assumptions {\bf (A2)}-{\bf (A3)} be satisfied and let assume that $q_j\equiv 0,\; \ell_j\equiv 1$ for all $j\in \calM$. If the underlying graph $\mathsf{G}$ is strongly connected, then $(\Sigma_{\mathsf{TN}})$ is boundary boundary approximately positive  controllable if and only if 
	\begin{equation}\label{Rank-condition}
		\overline{co}\left((\mathcal{I}^{out}_{\mathsf{w}})^{\top}\mathbb{A}^m KU_{+},\;\; m=0,1,\ldots,M-1\right)=\mathbb{R}^M_+,
	\end{equation}
	where $\mathbb{A}$ is the weighted transposed adjacency matrix of $\mathsf{G}$, see (\ref{weighted}).
\end{corollary}
{\it Proof}
Observe that, in this case, we have $\mathbb{A}(\mu)=e^{-\frac{\mu}{v}}\mathbb{A}$. Then, $$\Vert \mathbb{A}(\mu)\Vert\le e^{-\frac{\mu}{v}}\Vert \mathbb{A}\Vert< 1,$$ 
for all $\mu \ge 0$. Therefore, according to Theorem \ref{Theorem-transport}, $(\Sigma_{\mathsf{TN}})$ is boundary boundary approximately positive  controllable if and only if the following implication holds for all $\varphi\in L^q((0,1),\R^M)$:
$$
\langle (\mathcal{I}^{out}_{\mathsf{w}})^{\top}e^{-m\frac{\mu}{v}}\mathbb{A}^mKu,\varphi\rangle \leq 0, \qquad\, \forall\, u \in \R^n_+,\,m\in \N,\;\mu\ge 0\quad \Longrightarrow \quad \varphi \le  0,
$$
or, equivalently, if 
$$
\langle (\mathcal{I}^{out}_{\mathsf{w}})^{\top}\mathbb{A}^m Ku,\varphi\rangle \leq 0, \qquad\, \forall\, u \in \R^n_+,\,m\in \N,\;\mu\ge 0\quad \Longrightarrow \quad \varphi \le  0.
$$
On the other hand, by virtue of Cayley-Hamilton theorem, the power of $\mathbb{A}^m $ for all $m\ge M$ are linear combinations of the lower power, i.e., there exist real coefficients $a_m$, $m\in \{0,1,\ldots,M-1\}$ such that
$$
\mathbb{A}^M=\sum_{m=0}^{M-1}a_m\mathbb{A}^m.
$$
Moreover, the fact that the graph $\mathsf{G}$ is strongly connected yields that the matrix $\mathbb{A}$ is positive irreducible and hence the coefficients $a_m$ are nonegative. It follows that, the transport network systems $(\Sigma_{\mathsf{TN}})$ is boundary approximately positive  controllable if and only if the condition (\ref{Rank-condition}) holds.
\qed

\begin{example}\label{EX2}
	Here we consider a transport process on an $(N,N)$-directed cycle with a control acting on the starting node described as
	$$
	(\Sigma_{\mathsf{NP}})		
	\left\lbrace
	\begin{array}{lll}
		\frac{\partial}{\partial t} z(t,x)= c\frac{\partial }{\partial x}z(t,x),&t\geq 0, \;x\in (0,1),\cr  
		z(0,x)= f(x)\ge 0, &x\in (0,1),\;\,({\rm IC}) \cr
		z(t,1)= \mathbb{A}z(t,0)+Ku(t), & t\geq 0,\; \qquad({\rm BC})
	\end{array}
	\right.
	$$
	where  $c:={\rm diag}(c_j)_{j\in \calM}$ with $c_j>0$ for all $j\in \calM\in\{1,\ldots,N\}$, $K \in \calL(\R^N)$ with $(K)_{11}=b> 0$ and $(K)_{ij}=0$ otherwise, and $\mathbb{A}$ is the transposed adjacency matrix of the $(N,N)$-directed cycle.%, i.e., 
	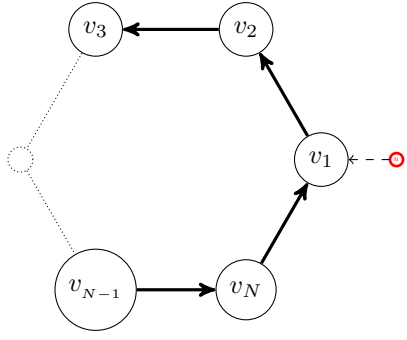
\begin{figure}[ht]%\label{Fig2}
		\begin{center}
			\begin{tikzpicture}
				\tikzstyle{every node}=[draw,shape=circle];
				%\node (v0) at (3:3) {$v_0$};
				\node[terminal] (P) at  (3,0) {$\textcolor{red}{\large u}$};
				\node (v1) at (0:2) {$v_1$};
				\node (v2) at (60:2) {$v_2$};
				\node (v3) at (2*60:2) {$v_3$};
				\node[densely dotted](v4) at (3*60:2) {};
				\node (v6) at (4*60:2) {$v_{_{N-1}}$};
				\node (v5) at (5*60:2) {$v_N$};
				\draw[dashed, ->] (P) -- (v1) ;
				\draw[fleche]  (v1) -- (v2);
				\draw[fleche] 	(v2) -- (v3);
				\draw[densely dotted] 	(v3) -- (v4);
				\draw[densely dotted]	(v4) -- (v6);
				%	\draw[densely dotted] 	(3*72:4) -- (3*72:2);
				\draw[fleche] 	(v6) -- (v5);
				\draw[fleche] 	(v5) -- (v1);
			\end{tikzpicture}
		\end{center}
		\caption{ The system of transport equations $(\Sigma_{\mathsf{NP}})$ on a cycle with a single input (represented by the read circle) acting on the starting node $v_1$.}
		\vspace{.1 in}
	\end{figure}
	
	Let $X=L^1([0,1])^N$, $\partial X=\R^{N}$, $U=\R$ and let us consider the operator $A_m:D(A_m)\to X$ defined by
	\begin{equation}\label{A_m-transport}
		A_{m}f:= c\partial_{x} f,\qquad
		D(A_{m}):= W^{1,1}([0,1])^N.
	\end{equation}	
	Moreover, we define the boundary operators $G$ and $\Gamma$ as follow 
	$$
	G f:=	f(1),\qquad
	\Gamma f :=\mathbb{A}f(0),\qquad f\in W^{1,1}([0,1])^N.
	$$
	It follows from Lemma \ref{Transport-wellposed}	that the system $(\Sigma_{\mathsf{NP}})$ is positively well-posed. In addition, we have that an $(N,N)$-directed cycle is strongly connected and
	$$
	\overline{co}\left(\mathbb{A}^m K \R_{+},\;\; m=0,1,\ldots,N-1\right)=\mathbb{R}^N_+.	
	$$
	Therefore, according to Corollary \ref{COR1}, the system $(\Sigma_{\mathsf{NP}})$ is boundary approximately positive  controllable.
\end{example}

\subsection{Heat conduction networks} %For a vertex v ∈ V, we denote by Γ (v) the set of all the edges in G incidentto v
We consider the following system of $M\ge 2$ coupled heat equations with controls in Robin boundary conditions,
\begin{align*}
	(\Sigma_{\mathsf{H}})		
	\begin{cases}
		\dfrac{\partial}{\partial t} z(t,x)= c\dfrac{\partial^2 }{\partial x^2}z(t,x)-q z(t,x),&t\geq 0, \;x\in (0,1),\cr  
		z(0,x)= h(x)\ge 0, &x\in (0,1),\;\,(\text{IC}) \cr
		\dfrac{\partial}{\partial x} z(t,1)= \mathbb{B}z(t,0)+Ku(t), & t\geq 0,\; \qquad\,(\text{BC})\\
		\dfrac{\partial}{\partial x} z(t,0)=0.
	\end{cases}
	%		\right.
\end{align*}
Here, $z_j(t,x)$ represents the heat distribution at time $t\ge 0$ and location $x\in (0,1)$ on an edge $\mathsf{e_j}$, where  $c:={\rm diag}(c_j)_{j\in \calM}$, $q:={\rm diag}(q_j)_{j\in \calM}$ with $c_j,\,q_j>0$ for all $j\in \calM:=\{1,\ldots,M\}$. The coupling between the $M$ heat equations is determined by the matrix $\B\in \calL(\R^M,\R^M)$, where we impose $K \in \calL(\R^n,\R^M)$ with $n\le M$ and $ 0\le u\in L_{loc}^2$
denotes the control matrix and the control vector, respectively. 

In order to apply the results of the previous sections, let $X=L^2([0,1])^M$, $\partial X=\R^M$, $U=\R^n$ and define the operator $A_m$ as
\begin{align}\label{A_m-heat}
	A_{m}h:= c\partial_{xx} h- q h,\quad
	D(A_{m}):= \{h\in W^{2,2}([0,1])^M:\; \partial_xh(0)=0\;\},
\end{align}	
where $c:={\rm diag}(c_j)_{j\in \calM}$ and $q:={\rm diag}(q_j)_{j\in \calM}$. Moreover, the boundary operators are given by 
\begin{align*}%\label{Operators-boundary}
	G h:=	\partial_xh(1),\qquad
	\Gamma h :=\mathbb{B}h(0),\qquad h\in W^{2,2}([0,1])^M.
\end{align*}

With the above notation, the system $(\Sigma_{\mathsf{H}})$ is rewritten in the form \eqref{S5.Sy1}.	

\begin{lemma}\label{generation}
	Let assume that $c_j,\,q_j>0$ for $j\in \calM$. Then the operator
	\begin{align}\label{operator-heat}
		A=A_m,\quad D(A):= \{h\in W^{2,2}([0,1])^M:\;\partial_xh(1)=0,\; \partial_xh(0)=0\},
	\end{align}
	generates a uniformly exponentially stable positive C$_0$-semigroup $\T$ on $X$ given by 
	\begin{align}\label{Expression-semigroup}
		(T(t)h)_j=\sum_{k=0}^{+\infty} e^{\lambda_{j,k} t} \langle h_j,\varphi_{k} \rangle \varphi_{k} ,
	\end{align}
	for all $j\in \calM$, where $ \lambda_{j,k} =-q_j-c_j k^2 \pi^2$ and $\varphi_{j,k}(x) =\cos( k\pi x)$ for $x\in [0,1]$. %\sqrt{c_j}
\end{lemma}
{\it Proof}
To prove our claim we shall	consider only one index $j\in \calM$ and in a similar way we deduce that of any $j\in \calM$. Indeed, let $j\in \calM $ be fixed and let $h_j\in D(A_j)$ such that $(Ah)_j=0$, then we have 
\begin{align*}
	\begin{cases}
		c_j	 \partial_{xx}h_j(x)-(q_j+\lambda_{j})h_j(x)=0,& x\in (0,1),	\\
		\dot{h}_j(1)=\dot{h}_j(0)=0.&
	\end{cases}
\end{align*}
So, solving the above eigenvalue problem, we get eigenvalues
\begin{align}\label{eigenvalues}
	\lambda_{j,k}=-q_j-c_jk^2\pi^2,\qquad k\in \N,
\end{align} 
and eigenfunctions $\varphi_{j,k}(x)= \cos( k\pi x)$ for $x\in (0,1)$. Then $(\varphi_{j,k})_{k\in \N}$ is an orthonormal basis in $L^2([0,1])$. Moreover, for every $g\in L^2([0,1])$ such that $(Ah)_j=g$, we have
\begin{align*}
	\begin{cases}
		c_j	\partial_{xx}h_j(x)-q_jh_j(x)=g(x),& x\in (0,1),	\\
		\dot{h}_j(1)=\dot{h}_j(0)=0.&
	\end{cases}
\end{align*}
which has the unique solution
\begin{align*}
	h_j(x)&=\frac{\int_{0}^{1}\sinh\left(\sqrt{\tfrac{q_j}{c_j}}(1-y)\right)c_j^{-1}g(y)dy}{-\sqrt{\tfrac{q_j}{c_j}} \sinh\left(\sqrt{\tfrac{q_j}{c_j}}\right)}\cosh(\sqrt{\tfrac{q_j}{c_j}})\\
	&+\left(\sqrt{\tfrac{q_j}{c_j}}\right)^{-1}\int_{0}^{x}\sinh\left(\sqrt{\tfrac{q_j}{c_j}}(x-y)\right)c_j^{-1}g(y)dy.
\end{align*}
Thus, $0\in \rho(A_j)$ and hence $A_j$ is diagonalizable, see \cite[Proposition 2.6.2]{TW}. Therefore, we get
\begin{align*}
	R(\mu,A_j)h_j=\sum_{k\in\N} \frac{1}{\mu-\lambda_{j,k}}\langle h_j,\varphi_{j,k} \rangle  \ge 0
\end{align*}
for all $h_j\in L^2_+([0,1])$ and ${\rm Re}\,\mu>-q_j$. Moreover, from \eqref{eigenvalues}, we have 
\begin{align*}
	\sup_k {\rm Re}\,\lambda_{j,k}\le -q_j, \qquad \forall\, j\in \calM.
\end{align*}
Hence, according to \cite[Proposition 2.6.5]{TW}, the operator $(A,\ker G)$ generates a positive C$_0$-semigroup $\mathbb{T}$ on $(L^2([0,1]))^M$ given by \eqref{Expression-semigroup} and $\omega_0(A)=\sup_j \sup_k {\rm Re}\,\lambda_{j,k}$. In addition, from \eqref{eigenvalues} we have $$\omega_0(A)\le - \underset{j\in \calM}{\sup }q_j<0.$$ Therefore, $\mathbb{T}$ is an uniformly exponentially stable positive semigroup. This completes the proof.
\qed

\begin{lemma}\label{Heat-wellposed}
	Let assume that $c_j,\,q_j>0$ for $j\in \calM$. Then the operator
	\begin{align}\label{CP2}
		\calA=A_{m},\qquad
		D(\calA)=\big\lbrace f\in D(A_m) :\; \partial_x f(1)=\mathbb{B}f(0)\;\big\rbrace,
	\end{align}
	generates a positive C$_0$-semigroup $\mathscr{T}:=(\mathscr{T}(t))_{t\geq 0}$ on $(L^2([0,1]))^M$. Moreover, 
	the system $(\Sigma_{\mathsf{H}})$ has a unique mild solution satisfying
	\begin{align}\label{variation-final}
		z(t)=\mathscr{T}(t)h+\int_{0}^{t} \mathscr{T}_{-1}(t-s)B K u(s) ds,\qquad t\ge 0,
	\end{align}
	for all $h\in X$ and $u\in L^2(\R_+,U)$. In addition, if $K$ is positive, then $z\in X_+$ for every initial data $h\in X_+$ and every input $u\in L_+^2(\R_+,U)$.
\end{lemma}
{\it Proof}
To prove that $\calA$ as defined in \eqref{CP2} generates a positive C$_0$-semigroup on $(L^2([0,1]))^M$, we shall use Theorem \ref{wellposedness}. In fact, it is clear that the operator $G$ ($:=\delta_1 \otimes \partial_{x}$) is surjective and by a simple computation we obtain that the Dirichlet operator $D_\mu$ associated to $G$ is given by 
\begin{align}\label{Dmu}
	(	D_\mu d)(x)&={\rm diag}\left(\xi_j(x)\right)_{j\in \calM}d, \qquad {\rm Re}\,\mu> -\underset{j\in \calM}{\sup }q_j:=\tilde{\mu}\\
	\xi_j(x)&:=\tfrac{\cosh\left(\sqrt{\tfrac{\mu+q_j}{c_j}} x\right)}{\sqrt{\tfrac{\mu+q_j}{c_j}}\sinh\left(\sqrt{\tfrac{\mu+q_j}{c_j}}\right)}.
\end{align}
for all $d\in \R^M$ and $x\in [0,1]$. Obviously, $D_\mu$ is positive for all ${\rm Re}\, \mu> \tilde{\mu}$. On the other hand, according to Lemma \ref{generation}, we have the operator $(A,\ker G)$ generates a positive C$_0$-semigroup $\mathbb{T}$ on $(L^2([0,1]))^M$ given by \eqref{Expression-semigroup}. So, in view of Theorem \ref{wellposedness}, it remains to show that $(A,B,\Gamma_{\vert \ker G})$ is a positive L$^2$-well-posed regular triplet on $(L^2([0,1]))^M,\R^M,\R^M$. In fact, let $C:=\Gamma_{\vert \ker G}$, let $\alpha>0$ and $0\le h\in D(A)$. Then,
\begin{align*}
	\int_{0}^{\alpha}\Vert CT(t)h\Vert_{\R^M}^2 dt&=\int_{0}^{\alpha}\Vert \mathbb{B}(T(t)h)(0)\Vert_{\R^M}^2 dt\\
	&\le \int_{0}^{\alpha}\sum_{j=1}^{N}\left\vert \sum_{k=0}^{+\infty} e^{\lambda_{j,k} t} \langle h_j,\varphi_{k} \rangle \right\vert^2 dt\\
	&\le \sum_{j=1}^{M} \sum_{k=0}^{+\infty}\int_{0}^{\alpha} e^{-2(q_j+c_j k^2 \pi^2)t} dt \sum_{k=0}^{+\infty} \left\vert\langle h_j,\varphi_{k} \rangle \right\vert^2\\
	&\le  \sum_{k=0}^{+\infty}  \sup_{j\in \calM}\frac{\alpha}{2(q_j+c_j k^2 \pi^2)}\sum_{j=1}^{M}\sum_{k=0}^{+\infty} \left\vert\langle h_j,\varphi_{k} \rangle \right\vert^2.
\end{align*}
Since $$\sum_{k=0}^{+\infty}  \frac{1}{2(q_j+c_j k^2 \pi^2)}<+\infty, \qquad \forall\, j\in \calM,$$
then there exists $\gamma(\alpha)>0$ such that 
\begin{align*}
	\int_{0}^{\alpha}\Vert CT(t)h\Vert_{\R^M}^2 dt\le \gamma(\alpha)\Vert h\Vert^2,
\end{align*}	
for all $0\le h\in D(A)(:=\ker G)$. Thus, $C$ is a positive L$^2$-admissible  observation operator for $A$. On the other hand, it is easy to see that the operator $A$  (see \eqref{operator-heat}) is self adjoint and the operator $B^*$, the adjoint operator of $B:=-A_{-1}D_0$, is given by 
\begin{align*}
	B^*h=h(1),\qquad \forall\, h\in (D(A))'.
\end{align*}
Similar argument to the admissiblity of $C$, it can be verified that $B^*$ is  positive L$^2$-admissible  observation operator for $A$. Hence, by duality (see e.g, \cite[Theorem 4.4.3]{TW}), we get that $B$ is a positive L$^2$-admissible control operator for $A$. The explicit expression of the operator $D_\mu$ (see \eqref{Dmu}) further yields that the transfer function of $(A,B,\Gamma_{\vert \ker G})$ is given by 
\begin{align}\label{transfer-fucntion}
	\mathbb{A}(\mu):=\Gamma D_\mu=\mathbb{B}{\rm diag}\left(\xi_j(0)\right)_{j\in \calM}, \qquad\forall\, {\rm Re}\, \mu >\tilde{\mu}.
\end{align}
Thus, $\mathbb{A}$ is analytic, bounded and positive for all ${\rm Re}\,\mu >\tilde{\mu}$. Moreover, we have $\lim_{{\rm Re}\mu\mapsto +\infty}\mathbb{A}(\mu) =0$, and hence according to Remark \ref{Regularity} $(A,B,\Gamma_{\vert \ker G})$ is a positive L$^2$-well-posed regular triplet on $(L^2([0,1]))^M,\R^M,\R^M$. It follows that there exists $\mu_0>0$ large enough such that $\Vert\mathbb{A}(\mu_0) \Vert<1$ and hence $r(\mathbb{A}(\mu_0))<1$. Then, in view of \cite[Lemma 4.1]{El} the identity $I_{\R^N}$ is a positive admissible feedback operator for $(A,B,\Gamma_{\vert \ker G})$. Therefore, according to Theorem \ref{wellposedness}, the operator $\calA$ from \eqref{CP2} generates a positive C$_0$-semigroup $\mathscr{T}:=(\mathscr{T}(t))_{t\geq 0}$ on $ X $ and the system of heat equations $(\Sigma_{\mathsf{H}})$ has a unique mild solution on $X$ given by the variation of constant formula \eqref{variation-final}, which is positive if $K$ is positive. %
\qed

\begin{lemma}\label{L-heat}
	Let assume that $c_j,\,q_j>0$ for $j\in \calM$. Then, the system of heat equations $(\Sigma_{\mathsf{H}})$ is approximately controllable with respect to positive controls if and only if 
	\begin{align}\label{F5}
		\bigcap\left\{\left(\mathbb{D}_\mu\mathbb{A}^m(\mu) K\R_{+}^n\right)^{\circ},\;\; m\in \N \right\}=\{0\}, \qquad \forall \,\mu\ge \mu_0,
	\end{align}
	for some $\mu_0>0$ large enough.
\end{lemma}
{\it Proof}
According to the proof of Lemma \ref{Heat-wellposed}, there exists $\mu_0>0$ large enough such that $r(\mathbb{A}(\mu_0))<1$. Thus, for all $\mu \ge \mu_0$ we have $r(\mathbb{A}(\mu))<1$, since the family $(\mathbb{A}(\mu))_{\mu>\mu_0}$ is positive and monotonically decreasing. By using Theorem \ref{Main-T1} together with the Neumann series representation of $(I-\mathbb{A}(\mu))^{-1}$ we obtain that the system $(\Sigma_{\mathsf{H}})$ is approximately controllable with respect to positive controls if and only if the condition \eqref{F5} holds for all $\mu \ge \mu_0$.
\qed

\begin{example}
	Here we consider a system of heat equations on a $(3,3)$-directed path (i.e., constituting of 3-directed edges connecting 3-vertices) with a control acting on the starting node described as
	\begin{align*}
		(\Sigma_{\mathsf{H}_{{\rm path}}})		
		\begin{cases}
			\dfrac{\partial}{\partial t} z(t,x)= c\dfrac{\partial^2 }{\partial x^2}z(t,x)-q z(t,x),&t\geq 0, \;x\in (0,1),\cr  
			z(0,x)= h(x)\ge 0, &x\in (0,1),\;\,(\text{IC}) \cr
			\dfrac{\partial}{\partial x} z(t,1)= \mathbb{B}z(t,0)+Ku(t), & t\geq 0,\; \qquad\,(\text{BC})\\
			\dfrac{\partial}{\partial x} z(t,0)=0,
		\end{cases}
		%		\right.
	\end{align*}
	where the coupling matrix $\B$ and the control matrix $K$ are given by, respectively, 
	\begin{align*}
		\mathbb{B}=\begin{pmatrix}
			0 & 0 & 0\\
			1 & 0 &  0\\
			0& 1 &  0
		\end{pmatrix}, \qquad K=\begin{pmatrix}
			b & 0 & 0  \\
			0 & 0 &  0 \\
			0& 0 &  0
		\end{pmatrix},\qquad b\in \R^*.
	\end{align*}

	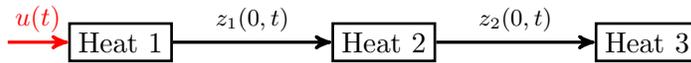
\begin{figure}[ht]\label{FigN}
		\begin{center}
			\begin{tikzpicture}[scale=0.5]
				%\node[terminal] (p) at (-19,0) {$\textcolor{red}{\textbf{+}}$};
				%	\node[element] (S) at (10,0) {\large Heat N};
				\node[element] (M) at (-15,0) {\large Heat 1};
				\node[element] (N) at (-8,0) {\large Heat 2};
				\node[element] (O) at (-1,0) {\large Heat 3};
				%		\node[element] (D) at (4,0) {\ldots};
				%\node[element] (F) at (4,-2) {$``\emph{v=y}``$};
				\draw[fleche] (M) -- (N) node[midway, above]{\small $z_1(0,t)$};
				\draw[fleche] (N) -- (O) node[midway, above]{\small $z_2(0,t)$};
				%	\draw[densely dotted] (O)-- (5,0) ;
				%	\draw[fleche] (5,0) --(S)node[midway, above]{\small $z_{N-1}(0,t)$};
				%	\draw[fleche] (F) -| (p);
				\draw[fleche1] (-18,0) -- (M) node[midway, above]{$\textcolor{red}{u(t)}$};
			\end{tikzpicture}
		\end{center}
		\caption{ The line graph representation of the system of heat equations $(\Sigma_{\mathsf{H}_{{\rm path}}})$.}
	\end{figure}
	
	We have the following negative controllability result.	
	\begin{lemma}\label{Lemma-path}
		Let assume that $c_j,\,q_j>0$ for $j\in \calM:=\{1,2,3\}$. Then, the system heat equations $(\Sigma_{\mathsf{H}_{{\rm path}}})$ is not approximately controllable with respect to a positive control acting on the starting node.
	\end{lemma}
	
	{\it Proof}
	In view of Lemma \ref{L-heat}, $(\Sigma_{\mathsf{H}_{{\rm path}}})$ is boundary approximately controllable with respect to positive controls if and only if 
	\begin{align}\label{path}
		\bigcap\left\{\left(\mathbb{D}_\mu\mathbb{A}^m(\mu) K\R_{+}\right)^{\circ},\;\; m\in \N \right\}=\{0\}, \qquad \forall \,\mu\ge \mu_0,
	\end{align}
	for some $\mu_0>0$ large enough, where
	\begin{align}\label{heat-Dlambda}
		\mathbb{A}(\mu):=\Gamma\mathbb{D}_\mu=\mathbb{B}{\rm diag}\left(\xi_j(0)\right)_{j=1,2,3},
	\end{align}
	for any $d\in \R^3$ , $x\in [0,1]$ and ${\rm Re}\,\mu>- \underset{j\in \calM}{\sup }q_j$. Let $\mu\ge \mu_0$, then by a direct computation we obtain
	\begin{align*}
		\mathbb{A}^2(\mu)= \begin{pmatrix}
			0 & 0 & 0\\
			0 & 0 &  0\\
			\xi_2(0)\xi_1(0) & 0 &  0
		\end{pmatrix},\qquad {\rm and }\qquad \mathbb{A}^m(\mu)=  0, \quad \forall\, m\ge 3.
	\end{align*}
	Therefore, the condition \eqref{path} is equivalent to the following implication:
	\begin{align}
		\langle \mathbb{H}u,h \rangle\le 0,\qquad \forall\, u\in \R_+\implies h=  0,
	\end{align} 
	where the operator $\mathbb{H}$ is given by 
	\begin{align}\label{H-operator}
		(\mathbb{H}u)(x)= \begin{pmatrix}
			\xi_1(x) bu\\
			\xi_2(x)\xi_1(0)bu \\
			\xi_3(x)\xi_2(0)\xi_1(0)bu
		\end{pmatrix}, \qquad \forall\, u\in \R_+,\,x\in (0,1).
	\end{align}
	Therefore, $(\Sigma_{\mathsf{H}_{{\rm path}}})$ is boundary approximately controllable with respect to positive controls if and only if, for all $\mu>\mu_0$,
	\begin{align}
		\langle \xi_1(\cdot) b ,\varphi_{1,k} \rangle< 0,\;\; \langle \xi_2(\cdot)\xi_1(0) b ,\varphi_{2,k} \rangle< 0,\; {\rm and }\;  \langle \xi_3(\cdot) \xi_2(0)\xi_1(0) b ,\varphi_{3,k} \rangle< 0,
	\end{align} 
	for $k\in \N$, where $\varphi_{j,k}$ (for $j=1,2,3$) is a Riesz basis in $(L^2(0,1))^3$ (see Lemma \ref{generation}). Moreover, a simple computation using the explicit expression of the function  $\xi_1(\cdot) $ yields  
	\begin{align}
		\langle \xi_1(\cdot) b ,\varphi_{1,k} \rangle	= \frac{(-1)^k}{\tfrac{\mu+q_j}{c_j}+(k\pi)^2}b,\qquad  k\in \N.
	\end{align}
	Therefore,  
	$$\langle \xi_1(\cdot) b ,\varphi_{1,k} \rangle< 0, \quad \forall\,k\in \N,\,\mu>\mu_0,$$%\Longrightarrow b=0 %\qquad  {\rm if\; and\; only\; if\;} \qquad b> 0 \quad {\rm and }\quad b< 0.
	which is impossible as $b\in \R^*$. This indicates that the heat equation on a path-like network is not approximately controllable with respect to a positive control acting on the starting node.
	\qed
	
	However, for positive input matrices, the following result shows that the system of heat equations $(\Sigma_{\mathsf{H}_{{\rm path}}})$ is boundary approximately positive  controllable.
	\begin{lemma}\label{Lemma-path1}
		Let assume that $c_j,\,q_j>0$ and $b>0$ for $j\in \calM$. Then the system $(\Sigma_{\mathsf{H}_{{\rm path}}})$ is boundary approximately positive  controllable.
	\end{lemma}
	{\it Proof}
	Let $\mu \ge \mu_0$ with $\mu_0>0$ large enough as in Lemma \ref{Lemma-path}. Then, according to Theorem \ref{Main-T2} and using the Neumann series representation of $(I-\mathbb{A}(\mu))^{-1}$ together with the computations from the proof of Lemma \ref{Lemma-path}, we get that the system $(\Sigma_{\mathsf{H}})$ is approximate positive controllable if and only if 
	\begin{align}\label{path1}
		\langle \mathbb{H}u,h \rangle\le 0,\qquad \forall\, u\in \R_+\implies h\le  0,
		%\bigcap\left\{\left(\mathbb{D}_\mu\mathbb{A}^m(\mu) K\R_{+}^4\right)^{\circ},\;\; m\in \N \right\}=-(L_+^2(0,1))^3, \qquad \forall \,\mu\ge \mu_0.
	\end{align}
	where $\mathbb{H}$ is given by \eqref{H-operator}. 
	
	To prove the claim of the lemma, we will argue by contradiction. In fact, let us assume that $\langle \mathbb{H}u,h \rangle\le 0$ for all $ u\in \R_+$ and $h>0$. Then, according to the calculations in the proof of Lemma \ref{Lemma-path}, we obtain
	\begin{align*}
		\langle \xi_1(\cdot) b ,\varphi_{1,k} \rangle \langle h_1,\varphi_{1,k} \rangle \le 0, \quad \forall\, k\in \N \qquad {\rm and }\qquad \langle h_1,\varphi_{1,k} \rangle>0, \quad \forall\, k\in \N. 
	\end{align*} 
	Thus, for even $k\in \N$ we obtain that: %$\langle h_1,\varphi_{1,k} \rangle>0$ and
	\begin{align*}
		\frac{b}{\tfrac{\mu+q_j}{c_j}+(k\pi)^2} \langle h_1,\varphi_{1,k}\rangle=\langle \xi_1(\cdot) b ,\varphi_{1,k} \rangle \langle h_1,\varphi_{1,k} \rangle   \le 0.
	\end{align*}
	As $b>0$ we get $ \langle h_1,\varphi_{1,k} \rangle   \le 0$ which is a contradiction. Therefore, the implication \eqref{path1} holds and hence the system of heat equations $(\Sigma_{\mathsf{H}_{{\rm path}}})$ is approximately positive controllable. 
	\qed
\end{example}

\begin{remark}
	Note that the results of Lemmas \ref{Lemma-path} and \ref{Lemma-path1} seem to be natural. Indeed, in Lemma \ref{Lemma-path}, we want to achieve a dense subspace of the entire state space $(L^2(0,1))^3$ by means of a positive control input through the leading (first) heat equation so that the boundary temperature of the first heat is fed into the boundary heat flux of the next one and so on. Therefore, the latter controllability property seems impossible, because the solutions of the heat equation $(\Sigma_{\mathsf{H}_{{\rm path}}})$ remain positive for any positive initial heat distribution. On the other hand,
	Lemma \ref{Lemma-path1} yields that if we start with positive initial heat distributions, we can reach approximately all positive heat distributions of  $(\Sigma_{\mathsf{H}_{{\rm path}}})$ by a positive control input through the leading heat equation.
\end{remark}

	\section{Conclusions}
In this paper, we obtained some general boundary approximate controllability criteria of infinite-dimensional control systems under positivity constraints . These results unify and extend some existing results to more general cases. Moreover, it has been shown that from our general results, some computable formulas
can be derived easily. For instance, assuming that the underlying graph is strongly connected, we obtained a Kalman-type rank condition for the boundary approximate controllability under positivity constraints on the control and state of transportation networks. An example of a system of transport equations on a cycle-like network with a single positive input acting on the starting node is also given. Furthermore, we proved the positivity of a coupled heat equations with controls in Robin boundary conditions. In particular, for a system of heat equations on a path-like network with a control acting on the leading heat equation, we  established the approximate controllability under positivity constraints on the control and state. However, we showed the lack of controllability under unilateral control-constraint.

Despite the fact that the applications have been focused on transport and heat systems, we conjecture that other potential candidates to these applications are more general parabolic equations and systems, hyperbolic systems and delayed systems.

	\appendix\label{Appendix}  %This command ends the counting of sections.
	\section*{Appendix}
	Here, we provide an appendix on technical lemmas needed for the proof of Corollary \ref{COR1}.
	
	The following result introduce a modification of the so called Sz\'{a}sz-Mirakjan operator, cf. \cite{Altomare}.%positive operators to prove a density. 
	\begin{lemma}\label{density0}
		Let $f\in \mathcal{C}_b([0,+\infty)):=\{f\in\calC([0,+\infty)):\;\exists\, \alpha\ge 0,\, \delta\ge 0\;{\rm such \, that }\; \vert f(x)\vert \le\delta e^{\alpha x} \} $ and define
		\begin{align}\label{Mirakjan}
			\mathscr{M}_n(f;x):=e^{-n\varphi(x)}\sum_{k=0}^{\infty} f(\varphi^{-1}(\tfrac{k}{n}))\frac{n^k}{k!}(\varphi(x))^k,\qquad n\ge 1,\; x\ge 0,
		\end{align}
		where $\varphi(x)=\frac{1}{v}(1-x)$ for $x\in [0,1]$ and $\varphi(x)=0$ for $x\ge 1$. Then the operators $\mathscr{M}_n$ are linear positive and for every $f\in \mathcal{C}_b([0,+\infty))$ we have 
		\begin{align*}
			\lim_{n\to+\infty}\mathscr{M}_n f=f,\qquad { \rm uniformly\, on} \;[0,1].
		\end{align*}
	\end{lemma}	
	{\it Proof}
	To prove our claim we will use the Korovkin theorem, see e.g. \cite{Altomare}. To this end, we have to prove that the operators $\mathscr{M}_n$ preserve the functions  $1$, $\varphi(x)$, $(\varphi(x))^2$. Indeed, a simple computation shows 
	\begin{align*}
		\mathscr{M}_n(1;x)&=1,\qquad
		\mathscr{M}_n(\varphi(x);x)= \varphi(x),
	\end{align*}
	and 
	\begin{align*}
		\mathscr{M}_n((\varphi(x))^2;x)&= e^{-n\varphi(x)}\sum_{k=0}^{\infty} \frac{k^2}{n^2}\frac{n^k}{k!}(\varphi(x))^k\\
		&=(\varphi(x))^2 + \frac{1}{n}\varphi(x),
	\end{align*}
	where we have used the fact that $\varphi^2(\varphi^{-1}(x))=x^2$. Therefore, from \cite[Theorem 4.1]{Altomare}, we get that
	\begin{align*}
		\lim_{n\to+\infty}\mathscr{M}_n f=f
	\end{align*}
	uniformly on $[0,1]$.
	\qed
	
	The last lemma at the hand, one can derive the following density result.
	\begin{lemma}\label{density1}
		Let $p\in [1,\infty)$ and $v>0$ be fixed. Then we have 
		\begin{align}\label{D1}
			\overline{co}\left(e^{-\frac{n}{v}(1-.)},\qquad n\in \N\right)=L^{p}_+([0,1].
		\end{align}
		%where we recall that $(L^{p}_+([-r,0]))^{\circ}=\left\{\varphi\in L^q([-r,0]),\; \langle  f,\varphi \rangle,\;\forall\, f\in L^p_+([-r,0])\right\}$.
		%for $\frac{1}{p}+\frac{1}{q}=1$.
	\end{lemma}	
	{\it Proof}
	Let $p,q\in [1,\infty)$ with $\frac{1}{p}+\frac{1}{q}=1$ and let $g\in L^q([0,1])$ such that $
	\int_{0}^{1} e^{-\frac{n}{v}(1-x)}g(x)\, dx\, \ge 0$ for all $n\in \N$. Let $0\le f\in \mathcal{C}([0,1])$ and define the function
	\begin{align*}
		h(x):=\begin{cases}
			f(x),& x\in [0,1],\cr
			f(1),& t\ge 1.
		\end{cases}
	\end{align*}
	Then $0\le h\in \mathcal{C}_b([0,+\infty))$ and 
	\begin{align*}\
		\int_{0}^{1} (\mathscr{M}_n h)(x)g(x)\,dx \,\ge 0.
	\end{align*}
	The continuity of $f$, Lemma \ref{density0} and the dominated convergence theorem further yield 
	\begin{align*}
		\int_{0}^{1}  f(x) g(x)\,dx\,\ge 0, \qquad \forall\, 0\le f\in \mathcal{C}([0,1]). 
	\end{align*}
	Moreover, since the positive cone in $\mathcal{C}([0,1])$ is dense in $L^p_+([0,1])$, we get that 
	\begin{align*}
		\int_{0}^{1}  f(x) g(x)\,dx\,\ge 0, \qquad \forall  f\in L^p_+([0,1]),
	\end{align*}
	and hence $g\ge 0$. 
	\qed
	%References

\end{document}